\documentclass[10pt, twoside]{article}
\usepackage{amssymb,amsmath}
\usepackage[amsmath,thmmarks,thref]{ntheorem}
\usepackage{graphics}
\usepackage{enumerate}
\usepackage{mathrsfs}
\title{Topological Representation of Geometric Theories}
\author{Henrik Forssell\thanks{Dept.\ of Informatics, University of Oslo, PO Box 1080 Blindern, N-0316 Oslo,
Norway. Email: jonf@ifi.uio.no}
}

\input{diagxy}
\xyoption{curve}
\newbox\anglebox 
\setbox\anglebox=\hbox{\xy \POS(75,0)\ar@{-} (0,0) \ar@{-} (75,75)\endxy}
 \def\pbangle{\copy\anglebox}
\newbox\angleboxr 
\setbox\angleboxr=\hbox{\xy \POS(0,0)\ar@{-} (0,75) \ar@{-} (75,0)\endxy}
 
\newbox\sanglebox 
\setbox\sanglebox=\hbox{\xy \POS(50,0)\ar@{-} (0,0) \ar@{-} (50,50)\endxy}
 
\newbox\sangleboxr 
\setbox\sangleboxr=\hbox{\xy \POS(0,0)\ar@{-} (0,50) \ar@{-} (50,0)\endxy}
 
\newbox\sangleboxf 
\setbox\sangleboxf=\hbox{\xy \POS(50,50)\ar@{-} (50,0) \ar@{-} (0,50)\endxy}
 
\newbox\angleboxf 
\setbox\angleboxf=\hbox{\xy \POS(75,75)\ar@{-} (75,0) \ar@{-} (0,75)\endxy}
 
\newbox\sangleboxfr 
\setbox\sangleboxfr=\hbox{\xy \POS(0,50)\ar@{-} (50,50) \ar@{-} (0,0)\endxy}
 
\newbox\angleboxfr 
\setbox\angleboxfr=\hbox{\xy \POS(0,75)\ar@{-} (75,75) \ar@{-} (0,0)\endxy}

\newdir{|>}{!/4.7pt/\dir{|}
        *:(1,-.2)\dir^{>}
        *:(1,+.2)\dir_{>}}

\def\embedd{\to/^{ (}->/}
\newcommand{\pair}[1]{\ensuremath{\langle {#1} \rangle}}
\newcommand{\homset}[3]{\ensuremath{\operatorname{Hom}_{#1}\!\left({#2},{#3}\right)}}

\newcommand{\mb}{\ensuremath{\mathbin}}

\newcommand{\cat}[1]{\ensuremath{\mathcal{#1}}}
\newcommand{\thry}[1]{\ensuremath{\mathbb{#1}}}
\newcommand{\synt}[2]{\ensuremath{\mathcal{#1}_{\mathbb{#2}}}}
\newcommand{\alg}[1]{\ensuremath{\mathbf{#1}}}
\newcommand{\mng}[1]{\ensuremath{\mathrm{#1}}}

\newcommand{\modin}[2]{\ensuremath{\mathrm{Mod}_{#1}(#2)}}
\newcommand{\topo}[1]{\ensuremath{\mathscr{#1}}}
\newcommand{\Sets}{\ensuremath{\mathbf{Set}}}

\newcommand{\Sh}[1]{\protect\ensuremath{\operatorname{Sh}\!\left(#1\right)}}

\newcommand{\Form}{\protect\ensuremath{\operatorname{Form}}}
\newcommand{\Mod}{\protect\ensuremath{\operatorname{Mod}}}

\newcommand{\sem}[1]{\ensuremath{[\![{#1}]\!]}}
\newcommand{\csem}[2]{\ensuremath{[\![{#1}\mb|{#2}]\!]}}

\newcommand{\fins}[1]{\exists {#1}\mathpunct .}

\newcommand{\theory}{\ensuremath{\mathbb{T}}}

\newcommand{\cterm}[2]{\ensuremath{\left \{ {#1}\ \; \vrule \; \ {#2}\right \}}}
\newcommand{\syntob}[2]{\ensuremath{[{#1}\;|\;{#2}]}}
\newcommand{\classtop}{\ensuremath{\mathbf{Set}[\theory]}}

\newcommand{\Eqsheav}[2]{\protect\ensuremath{\operatorname{Sh}_{#1}(#2)}}
\newcommand{\sox}[2]{\ensuremath{\csem{#1}{#2}_{M_{\theory}}}}
\newcommand{\bopen}[1]{\ensuremath{\langle\!\! \langle {#1} \rangle\!\! \rangle}}
\makeindex 
\begin{document}
%
%
\maketitle
\begin{abstract}
Using Butz and Moerdijk's topological groupoid representation of a topos with enough points, a `syntax-semantics' duality for geometric theories is constructed. The emphasis is on a logical presentation, starting with a description of the semantical topological groupoid of models and isomorphisms of a theory and a direct proof that this groupoid represents its classifying topos. Using this representation, a contravariant adjunction is constructed between theories and topological groupoids. The restriction of this adjunction yields a contravariant equivalence between theories with enough models and semantical groupoids. Technically a variant of the syntax-semantics duality constructed in \cite{awodeyforssell} for first-order logic, the construction here works for arbitrary geometric theories and uses a slice construction on the side of groupoids---reflecting the use of `indexed' models in the representation theorem---which in several respects simplifies the construction and allows for an intrinsic characterization of the semantic side.


\noindent \textbf{AMS classification codes:} 03G30; 18B25; 18C10; 18C50; 18F20.
\end{abstract}
\tableofcontents
\section{Introduction}
\label{introduction}

Grothendieck toposes can be regarded from both an `algebraic-geometric' point of view and, logically, from a `syntactical-semantic' point of view. If \theory\ is a geometric theory, i.e.\ a deductively closed set of sequents consisting of formulas constructed with the connectives $\top$, $\bot$, $\wedge$, $\exists$, and $\bigvee$ (where the latter is infinitary disjunction, but with the restriction that  no formula may contain more than finitely many free variables), then there exists a classifying topos, \classtop, with the property that the category of \theory-models in a topos, \cat{E}, is equivalent to the category of geometric morphisms from \cat{E} to \classtop\ (see \cite{elephant1}). There is a canonical `syntactical' construction of \classtop\ as sheaves on a site defined in terms of \theory.  On the other hand, for any topos \cat{E} one can define a geometric theory such that \cat{E} is the classifying topos of that theory. There is, therefore, a close connection between toposes and geometric theories. On the other hand, there is also a close connection between toposes and localic or topological groupoids; the category of equivariant sheaves on a localic groupoid is a topos, and any topos can be regarded as the topos of equivariant sheaves on some localic groupoid (see \cite{joyal:84}).  Moreover,  if the topos has enough points, then it can be regarded as the topos of equivariant sheaves on a topological groupoid, that is on a groupoid object in the category of topological spaces and continuous maps (see \cite{butz:98b}). Most toposes of interest have enough points, and we restrict attention here to those that do, and thus to topological groupoids (instead of localic ones). In logical terms, a topos regarded as a classifying topos \classtop\ of a theory \theory\ has enough points if and only if the theory has \emph{enough models},  in the sense that a sequent is provable in the theory if it is true in all models in the category of sets and functions, $\Sets$.

For any topos \cat{E} with enough points, then, one can regard \cat{E} as both the classifying topos of some theory (with enough models) and as a topos of equivariant sheaves for some topological groupoid, $\topo{G}=(G_1\rightrightarrows G_0)$
\[\classtop\simeq\cat{E}\simeq\Eqsheav{G_1}{G_0}\]
Now, the elements of the space of objects $G_0$ induce points of \Eqsheav{G_1}{G_0} and the elements of the space of arrows $G_1$ induce geometric transformations between them, and by the universal property of \classtop, these correspond to models and isomorphisms of \theory. Thus $\topo{G}$ can be regarded as a topological groupoid of \theory-models and isomorphisms (and the set of models $G_0$ is already `enough' for \theory, in the sense above). It is in this sense that  toposes can be regarded, `logically', from both a `syntactical' point of view and a `semantic' point of view.

Given a theory \theory, one can construct the classifying topos \classtop\ and then apply the Butz-Moerdijk construction from \cite{butz:98b} to obtain a topological groupoid, $G_1\rightrightarrows G_0$, consisting of models and isomorphisms of \theory. One can determine this groupoid directly from \theory, thus obtaining a direct semantical way to construct \classtop\ from \theory\ via its models. (At least in so far as the space of models is concerned, a related construction can also be found in the seminal work of Makkai and Reyes \cite{makkaireyes} who attribute the topology to \cite{hakim:72}). We state a (slightly simplified) variant of the construction in some detail here, and  supply a direct proof, using a minimum of topos theoretic machinery, that the topos of equivariant sheaves on it is the classifying topos of the theory. In the process we also obtain a direct description of the universal \theory-model, and relate it to Moerdijk's site description for equivariant sheaf toposes (see \cite{moerdijk:88}). It is clear that Butz and Moerdijk's representation result allows one to pass back and forth between  theories and their groupoid of models. But the precise formulation of the relationship between the category of theories and the category of groupoids involves making a number of choices, and no canonical or `best' formulation seems to have been given. In particular, one is faced with the question of which morphisms of groupoids to consider and how to characterize them, as well as how to characterize the groupoids that are in the image of the functor that sends a theory to its groupoid of models.  Also, one arbitrary choice occurs already  at the stage of constructing a topological groupoid of models for a given theory, involving a space, and thus a \emph{set} of models. The construction fixes a sufficiently---for that theory---large index set (or set of fresh constants, if you prefer) and considers models indexed by this set (term models). Thus one obtains a `semantical' functor that sends a theory to its groupoid of indexed models. For this reason, it is suggested here that the semantical functor should instead be seen as a functor into topological groupoids sliced over a particular groupoid `of sets' constructed from the index set. `Reindexing' functions between different choices of index sets then induce morphisms between their respectively induced groupoids, along which the respective slice categories can be compared.
For a fixed index set, with induced groupoid \topo{S} of sets, we construct  a contravariant adjunction between theories and topological groupoids sliced over \topo{S}. No (further) restriction of morphisms between groupoids is needed, and we can intrinsically characterize a full subcategory of `semantical' groupoids over \topo{S} such that the adjunction, when restricted to theories with enough indexed models and this subcategory, has isomorphisms as counit components and Morita equivalences as unit components. Thus we obtain a `syntax-semantics duality' between a category of theories and a slice category of groupoids of indexed models.
We  point out that we do not exploit the 2-categorical structure of theories and groupoids, preferring, for a simpler and more concrete presentation, to fix certain choices `on the nose' instead.

Section \ref{Section: Topological groupoids and categories of models} describes the semantic groupoid of a geometric theory and Section \ref{Section: Semantic representation of geometric theories} presents a simple and direct proof that taking equivariant sheaves on this groupoid results in the classifying topos of the theory. Section \ref{Section: Syntax-semantics duality} presents the adjunction between theories and the category of groupoids, and gives an intrinsic characterization of a full subcategory of groupoids to which is restricts to a duality in the sense above. A further restriction to coherent theories is also given.
%
%
The results and constructions of this paper resemble those in \cite{awodeyforssell}, in which a similar syntax-semantics duality was constructed for coherent theories which are decidable in the sense of having an inequality predicate. That restriction allowed for a somewhat simpler space of models, whereas here a topology (even) closer to the original in \cite{butz:98b} is used. The main difference from \cite{awodeyforssell} is, otherwise, in the use of a slice category of groupoids in order to obtain a much simpler set up for the duality and a characterization of groupoids of (indexed) models. (There is also a difference in perspective; in \cite{awodeyforssell} the main purpose is to show how the Butz-Moerdijk representation result can be reformulated and used to generalize the classical Stone duality from propositional to first-order logic. The groupoid \topo{S} of sets occurs there in the role of a dualizing object. Here, the purpose is to show the advantages simply slicing over \topo{S}.) Both \cite{awodeyforssell} and the current paper are based on the author's PhD dissertation \cite{phd}.


%
\section{Topological groupoids of models}
\label{Section: Topological groupoids and categories of models}
\subsection{The logical topology}
\label{Subsection: The logical topology}
Let a  signature (first-order with equality) $\Sigma$ be given. We assume for simplicity that $\Sigma$ is single-sorted (with no loss of generality for our purposes in so far as any geometric theory is Morita equivalent to a geometric theory with only one sort, see \cite[D1.4.13]{elephant1}).  Choose an infinite cardinal, $\kappa\geq |\Sigma|+\aleph_0$. Fix an  `index' set \thry{S} of size at least $\kappa$. For instance one can set $\thry{S}$ to be a sufficiently large initial segment of the set theoretical universe. Let $M_{\Sigma}$ be the set of $\Sigma$-structures the underlying sets of which are quotients of subsets of \thry{S}. We refer to such structures as (\thry{S}-)\emph{indexed}. Recall that a geometric theory is said to have \emph{enough models} if every geometric sequent which is true in all models is provable in \theory.
\begin{Lemma}\label{Lemma: Enumerated models reflect covers}
For any geometric theory, \theory, over $\Sigma$, if \theory\ has enough models then it has enough indexed models, i.e.\ models in $M_{\Sigma}$.
\begin{Proof}By the downward L\"{o}wenheim-Skolem theorem. Any elementary embedding is pure (in the sense of \cite{adamekandrosicky:94}) and so preserves and reflects geometric formulas and reflects the truth of geometric sequents. Given a \theory-model in which a geometric sequent, $\sigma$, is false, choose a witness of the falsehood and a smaller than $\kappa$ elementary substructure around it. This is then isomorphic to an indexed \theory-model in which $\sigma$ is false.
\end{Proof}
\end{Lemma}
We begin by defining a topology on $M_{\Sigma}$, as well as on the set of $\Sigma$-structure isomorphisms between them. We use boldface to indicate ordered finite lists or tuples, and write $\star$ for the empty list. We also use boldface to indicate structures and models, \alg{M},\alg{N}, and structure isomorphisms, \alg{f},\alg{g}. For a structure, \alg{M}, in $M_{\Sigma}$, the elements of \alg{M} are equivalence classes of elements of \thry{S}, which we write as $[a]$ (and $[\alg{a}]$ for a list) with (implicit) reference to the equivalence relation on the underlying set of \alg{M}.
\begin{Definition}\label{Definition> Logical topology}
The \emph{logical topology} on $M_{\Sigma}$ is the coarsest topology containing all sets of the following form:
\begin{enumerate}
\item For each element, $a\in \thry{S}$, the set \[\bopen{a}:=\cterm{\alg{M}\in M_{\Sigma}}{[a]\in {\alg{M}}}\]
(i.e.\ structures, \alg{M}, such that the underlying set of \alg{M} is a quotient of a subset of \thry{S} which contains $a$).
\item For each $n$-ary relation symbol, $R$, and $n$-tuple, $\alg{a}$, of elements of \thry{S}, the set \[\bopen{R,\alg{a}}:=\cterm{\alg{M}\in M_{\Sigma}}{[\alg{a}]\in \sem{R}^{\alg{M}}}\]
    (This includes equality and nullary relation symbols; we treat the extension of a sentence in a model as a subset of the distinguished terminal object, so if $R$ is a nullary relation symbol, then  $\bopen{R, \star}=\cterm{\alg{M}\in M_{\Sigma}}{\alg{M}\models R}$.)
\item For each $n$-ary function symbol, $f$, and list of n+1 elements, $\alg{a}, b\in \thry{S}$ the set $\bopen{f(\alg{a})=b}:=\cterm{\alg{M}\in M_{\Sigma}}{\sem{f}^{\alg{M}}([\alg{a}])=[b]}$.
\end{enumerate}
Let $I_{\Sigma}$ be the set of isomorphisms between the structures in $M_{\Sigma}$, with \emph{domain} and \emph{codomain} functions $d,c: I_{\Sigma}\rightrightarrows M_{\Sigma}$. The \emph{logical topology} on $I_{\Sigma}$ is the coarsest such that both $d$ and $c$ are continuous and containing all sets of the following form:
\begin{enumerate}[(i)]
\item For each pair of elements, $\pair{a,b}$, of \thry{S} the set
 \[\bopen{a\mapsto b}:=\cterm{\alg{f}\in I_{\Sigma}}{[a]\in d(\alg{f})\wedge [b]\in c(\alg{f})\wedge \alg{f}([a])=[b]}\]
\end{enumerate}
\end{Definition}
It is straightforward to verify (see Lemma \ref{Lemma: H continuous category}) that the composition map ($m$); the mapping ($e$) of a structure to its identity morphism ; and the mapping ($i$) of a isomorphism to its inverse  are continuous functions. Accordingly, we have a topological groupoid which we call the \emph{topological groupoid of models and isomorphisms}, $\topo{M}_{\Sigma}$:
%
%
%
%
%
%
%
%
%
%
\[\bfig
%
%
\morphism<750,0>[I_{\Sigma}\times_{M_{\Sigma}} I_{\Sigma}`I_{\Sigma};m]
%
%
\morphism(750,0)|a|/@{>}@<5pt>/<750,0>[I_{\Sigma}`M_{\Sigma};d]
\morphism(750,0)|m|/@{<-}/<750,0>[I_{\Sigma}`M_{\Sigma};e]
\morphism(750,0)|b|/@{>}@<-5pt>/<750,0>[I_{\Sigma}`M_{\Sigma};c]
\Loop(750,0)I_{\Sigma}(ur,ul)_{\mng{i}}
\efig\]
We are interested in the full topological subgroupoids formed by models of theories  over $\Sigma$. We consider a \theory-model, \alg{M}, simultaneously as a \cat{L}-structure satisfying \theory\ and as a functor $\alg{M}:\synt{C}{T}\to \Sets$, where \synt{C}{T} is the small syntactical category of \theory\ (see Section \ref{Subsection: Theories and Categories}). Accordingly, we consider a homomorphism between models simultaneously as a natural transformation of functors.
\begin{Remark}\label{Remark: Quotient str are really congruence str}
A structure in $M_{\Sigma}$ can be thought of in three ways. We can think of it, as it is defined, as a structure with underlying set a quotient $A/_{\sim}$ of a subset $A\subseteq \thry{S}$. But it can also, of course, be thought of as a model which has underlying set $A\subseteq \thry{S}$ and where equality is interpreted as the equivalence relation $\sim$. Isomorphisms between such models would then be appropriate relations. Under the latter point of view, the logical topology is a topology of `finite information' about the structure, in the sense that a basic open will specify a finite list of elements of the structure as well as finitely many conditions (in terms of $\Sigma$) that they satisfy (this is more explicitly the perspective in \cite{awodeyforssell}, where structures are subsets of the index set and not quotients of such subsets). Finally, one can consider the index set \thry{S} as a set of new constants and  a structure in $M_{\Sigma}$  as a form of term structure. The latter is the perspective on the space of models constructed in chapter 6 of \cite{makkaireyes}. The topology in \cite{butz:98b} differs mainly in requiring that the equivalence classes should be infinite.
\end{Remark}
%
%
\subsection{Theories and models}
\label{Subsection: Theories and models}

Let \theory\ be a   geometric theory over the signature $\Sigma$.  Set $M_{\theory}\subseteq M_{\Sigma}$ to be the set of \theory-models in $M_{\Sigma}$. Set $I_{\theory}\subseteq I_{\Sigma}$ to be the set of isomorphisms between models in $M_{\theory}$, and denote the resulting groupoid by $\topo{M}_{\theory}$.
We equip the sets $M_{\theory}\subseteq M_{\Sigma}$ and $I_{\theory}\subseteq I_{\Sigma}$ with the subspace topologies from $M_{\Sigma}$ and $I_{\Sigma}$, respectively. Recall that the logical topology is determined by the signature $\Sigma$, not the logical formulas over $\Sigma$. It is, however, convenient to note that a \emph{basic} open set of $M_{\Sigma}$ or $M_{\theory}$ can be presented in the form
\begin{equation}\label{Equation: Basic open set}\bopen{\syntob{\alg{x}}{\phi},\alg{a}}=\cterm{\alg{M}\in M_{\theory}}{\alg{a}\in \csem{\alg{x}}{\phi}^{\alg{M}}}\end{equation}
where \syntob{\alg{x}}{\phi} is a Horn formula and $\alg{a}\in \thry{S}$. A straightforward induction on formulas shows that for any geometric formula, $\phi$, with free variables in $\alg{x}$, the set defined by (\ref{Equation: Basic open set}) is open. We write this out for reference.
\begin{Lemma}\label{Lemma: The logical topologies are all the same}
Sets of the form
\[\bopen{\syntob{\alg{x}}{\phi},\alg{a}}=\cterm{\alg{M}\in M_{\theory}}{\alg{a}\in \csem{\alg{x}}{\phi}^{\alg{M}}}\]
form a basis for the logical topology on $M_{\theory}$, with $\phi$ ranging over all Horn formulas over $\Sigma$ or, redundantly, over all Cartesian, regular, coherent or geometric formulas.
\end{Lemma}
Similarly, a basic open set of $I_{\theory}$ can be presented in the following form:
\\
\begin{equation}\label{Equation: V array}
\left(\begin{array}{c}
\syntob{\alg{x}}{\phi},\alg{a}  \\
\alg{b} \mapsto \alg{c}   \\
 \syntob{\alg{y}}{\psi},\alg{d}
\end{array}\right)
=d^{-1}(\bopen{\syntob{\alg{x}}{\phi},\alg{a}})\bigcap\bopen{\alg{b} \mapsto \alg{c}}\bigcap c^{-1}(\bopen{\syntob{\alg{y}}{\psi},\alg{d}})
\end{equation}
where $\phi$ and $\psi$ are at least Horn formulas. We think of such a presentation of a basic open set as having a \emph{domain}, a \emph{preservation}, and a \emph{codomain condition.}
\begin{Lemma}\label{Lemma: H continuous category}
$\topo{M}_{\theory}$ is a topological groupoid, i.e.\ a groupoid object in the category \alg{Sp} of topological spaces and continuous functions:
%
%
%
%
%
%
%
%
%
%
\[\bfig
%
%
\morphism<750,0>[I_{\Sigma}\times_{M_{\Sigma}} I_{\Sigma}`I_{\Sigma};m]
%
%
\morphism(750,0)|a|/@{>}@<5pt>/<750,0>[I_{\Sigma}`M_{\Sigma};d]
\morphism(750,0)|m|/@{<-}/<750,0>[I_{\Sigma}`M_{\Sigma};e]
\morphism(750,0)|b|/@{>}@<-5pt>/<750,0>[I_{\Sigma}`M_{\Sigma};c]
\Loop(750,0)I_{\Sigma}(ur,ul)_{\mng{i}}
\efig\]
\begin{Proof}
The domain and codomain functions $d,c:I_{\theory}\two M_{\theory}$ are
continuous by definition. Checking that a subbasic open set $\bopen{a\mapsto b}\subseteq I_{\Sigma}$ is pulled back to an open set by the remaining inverse, identity, and composition maps,  we have $i^{-1}(\bopen{a\mapsto b})=\bopen{b\mapsto a}$ and $e^{-1}(\bopen{a\mapsto b})= \bopen{\syntob{x,y}{x=y},a,b}$, and, finally
\[m^{-1}(\bopen{a\mapsto b})=\bigcup_{c\in \thry{S}}\bopen{c\mapsto b}\times_{M_{\theory}}\bopen{a\mapsto c}. \]
\end{Proof}
\end{Lemma}
We refer to  $\topo{M}_{\theory}$ as the \emph{topological groupoid of models (and isomorphisms} of \theory. Before proceeding, we note that  the spaces $M_{\theory}$ and $I_{\theory}$ are sober. Recall, e.g.\ from \cite{johnstone:82}, that a completely prime filter in the frame of open sets of a topological space $X$ is a filter of opens, $\cat{F}\subseteq \cat{O}(X)$, with the property that if a union of a set of opens is in the filter, $\bigcup S\in\cat{F}$, then for some element, $U\in S$, we have that $U\in \cat{F}$. Recall that a space is sober if every completely prime filter of opens is the neighborhood filter of a unique point.
\begin{Proposition}\label{Proposition: Sober spaces of models}
The spaces $M_{\theory}$ and $I_{\theory}$ are sober.
\begin{Proof}
We prove that $M_{\theory}$ is sober, $I_{\theory}$ is similar.  Let a completely prime filter $F$ of open subsets of $M_{\theory}$ be given. Let $A\subseteq \thry{S}$  be the set
\[A:=\cterm{a\in \thry{S}}{\bopen{\syntob{x}{\top}, a}\in F}\]
Define an equivalence relation, $\sim$, on $A$ by
\[a\sim b\ \ \Leftrightarrow\ \ \bopen{\syntob{x,y}{x=y}, a,b}\in F\]
Define a $\Sigma$-structure on $A/_{\sim}$ by interpreting a relation symbol $R$ as the set
\[\sem{R}:=\cterm{[\alg{a}]}{\bopen{R,\alg{a}}\in F}\]
(note that $\bopen{R,\alg{a}}\in F$ and $\alg{a}\sim\alg{b}$ implies  $\bopen{R,\alg{b}}\in F$ since $F$ is closed under finite intersection) and a function symbol $f$ as the function
\[[\alg{a}]\mapsto [b]\ \ \Leftrightarrow\ \ \bopen{f(\alg{a})=b}\in F\]
We show that this  is well-defined. If $[\alg{a}]\in A/_{\sim}$ then $\bopen{\syntob{\alg{x}}{\top}, \alg{a}}\in F$, and
\[\bopen{\syntob{\alg{x}}{\top}, \alg{a}}\subseteq \bopen{\syntob{\alg{x}}{\fins{y}f(\alg{x})=y},\alg{a}}=\bigcup_{b\in \thry{S}}\bopen{\syntob{\alg{x},y}{f(\alg{x})=y},\alg{a},b}\]
whence there exists $b\in\thry{S}$ such that $\bopen{\syntob{\alg{x},y}{f(\alg{x})=y},\alg{a},b}=\bopen{f(\alg{a})=b}\in F$, so $f$ is total. If $\bopen{f(\alg{a})=b'}\in F$ then $\bopen{f(\alg{a})=b}\cap\bopen{f(\alg{a})=b'}\subseteq \bopen{\syntob{x,y}{x=y}, b,b'}$, so $f$ is functional.
We have defined, therefore, a structure $\alg{A}\in M_{\Sigma}$. By a straightforward induction on $\syntob{\alg{x}}{\phi}$, we have that
\[\alg{A}\vDash \phi([\alg{a}])\ \ \Leftrightarrow\ \ \bopen{\syntob{\alg{x}}{\phi}, \alg{a}}\in F\]
If \theory\ proves the sequent $\phi\vdash_{\alg{x}}\psi$ then $\bopen{\syntob{\alg{x}}{\phi},\alg{a}}\subseteq\bopen{\syntob{\alg{x}}{\psi},\alg{a}}$ for all $\alg{a}\in \thry{S}$, and so \alg{A} is a \theory-model. It is clear that $F$ is the neighborhood filter of \alg{A}, and that \alg{A} is unique with this property.
\end{Proof}
\end{Proposition}
Finally, we note that $\topo{M}_{\theory}$ is in fact an open topological groupoid, in the sense that the domain and codomain maps $d,c:I_{\theory}\rightrightarrows M_{\theory}$ are open maps. We need the following technical lemma (which will be much used also further down).
\begin{Lemma}\label{Lemma: Star lemma}
For any \alg{M} in $M_{\theory}$, any finite list $a_1,\ldots,a_k\in\thry{S}$ such that $[a_1],\ldots,[a_k]\in\alg{M}$, and any finite list of distinct elements $b_1,\ldots,b_k\in\thry{S}$, there exists a model \alg{N} in $M_{\theory}$ and an isomorphism $\alg{f}:\alg{M}\rightarrow\alg{N}$ such that
\[
\alg{f}([a_i])=[b_i]
\]
for all $1\leq i\leq k$.
\begin{Proof} Write the underlying set of \alg{M} as $A/_{\sim}$ with $A\subseteq \thry{S}$. Since
$b_1,\ldots,b_k$ are all distinct, $b_i\mapsto[a_i]$ defines a
partial function from \thry{S} to ${A}/_{\sim}$. Choose any
surjective  extension $p:\thry{S}\twoheadrightarrow
{A}/_{\sim}$, and write $\equiv$ for the induced equivalence
relation on \thry{S}. Then $[a]\mapsto p^{-1}{[a]}$
defines a bijection $f:A/_{\sim}\rightarrow
\thry{S}/_{\equiv}$ such that $f([a_i])=[b_i]$ for all
$1\leq i\leq k$. Thus we can let \alg{N} be the \theory-model
induced by $f$ and \alg{M} on $\thry{S}/_{\equiv}$.
\end{Proof}
\end{Lemma}
\begin{Proposition}\label{Proposition: I is an open groupoid}
The topological groupoid $\topo{M}_{\theory}$ of models and isomorphisms of \theory\ is an open topological groupoid.
\begin{Proof}
Suffice to show that, say, the domain map is open. Let
\[V= \left(\begin{array}{c}
\syntob{\alg{x}}{\phi},\alg{a}  \\
\alg{b} \mapsto \alg{c}   \\
 \syntob{\alg{y}}{\psi},\alg{d}
\end{array}\right) \subseteq I_{\theory}\]
be given. Assume, without loss of generality, that \alg{d} is a list of distinct elements. We can also assume that no element of \thry{S} occurs more
than once in the tuple $\alg{c}$ (since we can add identity statements to the domain condition to cut down to a single occurence). Let an isomorphism
$\alg{f}:\alg{M}\rightarrow \alg{N}$ in $V$ be given. Choose, for each element $d$ in the tuple $\alg{d}$, an element $k$ such that $\alg{f}([k])=[d]$ in such a way that if for some $c$ in the tuple $\alg{c}$ we have $d=c$ then
$k=b$. Now $\alg{M}\in \bopen{\syntob{\alg{x},\alg{y},\alg{z}}{\phi\wedge \psi},\alg{a},\alg{k},\alg{b}}$. And if $\alg{K}\in \bopen{\syntob{\alg{x},\alg{y},\alg{z}}{\phi\wedge \psi},\alg{a},\alg{k},\alg{b}}$, then we can find a \theory-model $\alg{L}$ and isomorphism $g:\alg{K}\rightarrow  \alg{L}$ such
that $g(\alg{[k]})=\alg{[d]}$ and $g(\alg{[b]})=\alg{[c]}$, by Lemma \ref{Lemma: Star lemma}. Hence $\alg{K}\in
d(V)$.
\end{Proof}
\end{Proposition}
\section{Semantic representation of geometric theories}
\label{Section: Semantic representation of geometric theories}
\subsection{Equivariant sheaves on topological groupoids}
\label{Subsection: Equivariant sheaves on topological groupoids}
We briefly recall the essentials concerning the topos of equivariant sheaves on a topological groupoid, see \cite{elephant1}, \cite{moerdijk:88}, and \cite{moerdijk:90} for more. Consider a topological groupoid \topo{G}:
\begin{equation}\bfig
%
%
\morphism<750,0>[G_1\times_{G_0} G_1`G_1;m]
%
%
\morphism(750,0)|a|/@{>}@<5pt>/<750,0>[G_1`G_0;d]
\morphism(750,0)|m|/@{<-}/<750,0>[G_1`G_0;e]
\morphism(750,0)|b|/@{>}@<-5pt>/<750,0>[G_1`G_0;c]
%
\place(750,100)[\curvearrowleft^i]
\efig\end{equation}
The objects of the category of \emph{equivariant sheaves}, \Eqsheav{G_1}{G_0}, on \topo{G} are pairs \pair{r:R\rightarrow G_0,\rho} where $r$ is a local homeomorphism---i.e.\ an object of \Sh{G_0}---and $\rho$ is a continuous action, i.e.\ a continuous map
\[\rho:G_{1}\times_{G_{0}}R\to R\]
(with the pullback being along the domain map) such that $r(\rho(f,x))=c(f)$, satisfying  unit and  composition axioms:
\[\rho(1_{r(x)},x)=x\ \ \ \ \ \  \rho(g\circ f,x)=\rho(g,\rho(f,x))\]
%
A morphism of equivariant sheaves is a morphism of sheaves commuting with the actions. Recall that the category, \Eqsheav{G_1}{G_0}, of equivariant sheaves on \topo{G} is a (Grothendieck) topos, and that the forgetful functors of forgetting the action,
$\Phi:\Eqsheav{G_1}{G_0}\to \Sh{G_0}$
and of forgetting the topology,
$\Psi:\Eqsheav{G_1}{G_0}\to \Sets^{\topo{G}}$
are both conservative inverse image functors. Explicitly, $\Psi$ sends an equivariant sheaf, \pair{r:R\rightarrow G_0,\rho}, to the functor which sends an arrow $f:x\rightarrow y$ in $G_1$ to the function $\rho(f,-):r^{-1}(x)\rightarrow r^{-1}(y)$.
%

A \emph{continuous functor}, or morphism of topological groupoids, $f:\topo{G}\to\topo{H}$, i.e.\ a morphism of groupoid objects in \alg{Sp}
\[\bfig
\square(0,0)|allb|/>`@{>}@<-3pt>`@{>}@<-3pt>`>/<600,300>[ G_{1}`H_{1}`G_{0}`H_{0};f_1`d`d`f_0]
\square(0,0)|arrb|/>`@{>}@<3pt>`@{>}@<3pt>`>/<600,300>[ G_{1}`H_{1}`G_{0}`H_{0};f_1`c`c`f_0]
\efig\]
induces a \emph{geometric morphism}, $f:\Eqsheav{G_1}{G_0}\to\Eqsheav{H_1}{H_0}$, that is, a pair of adjoint functors,
\[\bfig
\morphism|b|/{@{>}@/_1em/}/<1000,0>[\Eqsheav{G_1}{G_0}`\Eqsheav{H_1}{H_0};f_*]
\morphism/{@{<-}@/^1em/}/<1000,0>[\Eqsheav{G_1}{G_0}`\Eqsheav{H_1}{H_0};f^*]
\place(500,0)[\bot]
\efig\]
consisting of a \emph{direct image} functor $f_*$ and an \emph{inverse  image} functor $f^*$. The inverse image functor works by pullback in the expected way, and preserves finite limits (and therefore, being a left adjoint, geometric logic).
%
%
%
%
%

A site description for the topos of equivariant sheaves on an open localic groupoid is given by Moerdijk in \cite{moerdijk:88}, and we use it for the case  of equivariant sheaves on an open topological groupoid here (a detailed expose of the site description for this case can be found in \cite{forssell:subgpds}). Briefly,  let \topo{G} be an  open topological groupoid; let $N\subseteq G_1$ be an open subset of arrows that is closed under inverses and compositions; and let $U=d(N)=c(N)\subseteq G_0$. Then
\[d^{-1}(U)/_{\sim_N}\to^c G_0\]
is an equivariant sheaf over $G_0$, denoted \pair{\topo{G},U,N}, where $f\sim_N g$ iff $c(f)=c(g)$ and $g^{-1}\circ f \in N$. The action is defined by composition. The set of objects of this form is a generating set for \Eqsheav{G_1}{G_0}; briefly because if \pair{\rho, r:R\rightarrow G_0} is an equivariant sheaf and $s:U\rightarrow R$ is a continuous section,
\[\bfig
\dtriangle/<-`>`/<350,350>[R`U`G_0;s`r`]
\place(175,0)[\subseteq]
\efig\]
then $N=\cterm{f\in d^{-1}(U)\cap c^{-1}(U)}{\rho(f,s(d(f))=s(c(f))}$ is an open set of arrows closed under inverses and compositions. The map $e:U\rightarrow d^{-1}(U)$ defined by $x\mapsto [1_x]$ is a continuous section,  and $s$ lifts to a morphism $\hat{s}:\pair{\thry{G},U,N}\rightarrow \pair{\rho, r:R\rightarrow G_0}$ such that $s=\hat{s}\circ e$.  Refer to the full subcategory of  objects of the form \pair{\thry{G},U,N}  as the \emph{Moerdijk site} for \Eqsheav{G_1}{G_0} (the implicit coverage is the canonical one inherited from \Eqsheav{G_1}{G_0}). The following properties of Moerdijk sites will be of use and we state them in a single lemma here for reference (proofs can be found in the given references).
\begin{Lemma}\label{Lemma: GUN subobjects}Let \topo{G} be an  open topological groupoid. The Moerdijk site of \Eqsheav{G_1}{G_0} is closed under subobjects. In particular, let $N\subseteq G_1$ be an open subset of arrows that is closed under inverses and compositions, and let $U=d(N)=c(N)\subseteq G_0$. Then the frame of subobjects of the object \pair{\topo{G},U,N} in \Eqsheav{G_1}{G_0} is isomorphic to the frame of open subsets of $U$ that are closed under $N$, with such a $V\subseteq U$ corresponding to the subobject $d^{-1}(V)/_{\sim_{N\upharpoonright_V}}$.
\end{Lemma}

\subsection{Equivariant sheaves on the space of models}
\label{Section: Sheaves on the space of models}
Fix a geometric theory \theory\ with enough \thry{S}-indexed models and let $\topo{M}_{\theory}$ be its topological groupoid of models and isomorphisms, as in Section \ref{Subsection: The logical topology}. \theory\ has a classifying topos, \classtop, with the (defining) universal property that for any topos, \cat{E}, the category of \theory-models in \cat{E} is equivalent to the category of geometric morphisms from \cat{E} to \classtop,
\[\modin{\theory}{\cat{E}}\simeq\homset{\mathfrak{Top}}{\cat{E}}{\classtop}\]
and from which \theory\ can be recovered up to Morita equivalence.  The current section presents a direct---and alternative to that which can be found in \cite{butz:98b}---proof  that  the topos of equivariant sheaves on $\topo{M}_{\theory}$ is, in fact, (equivalent to) the classifying topos of \theory, thus yielding a semantic groupoid representation of \classtop\ (supplementing the standard syntactical construction.) In the process, we obtain a concrete description of the universal \theory-model in \Eqsheav{I_{\theory}}{M_{\theory}}, and it is shown that this model is a dense subcategory of the Moerdijk site of \Eqsheav{I_{\theory}}{M_{\theory}}. The proof presented here follows three steps: From the fact that \theory\ has enough models in $M_{\theory}$ it follows that there is a conservative embedding of $\cat{M}_d:\synt{C}{T}\to\Sets/_{M_{\theory}}$ and thus a surjective geometric morphism $m_{d}:\Sets/_{M_{\theory}}\epi\classtop$. The functor $\cat{M}_d$ is factored, first, through the category of sheaves on $M_{\theory}$
(equipped with the logical topology) and, second, through the category of equivariant sheaves on \Eqsheav{I_{\theory}}{M_{\theory}}:
\[\bfig
\qtriangle/>`>`<-/<750,300>[\synt{C}{T}`\Sh{M_{\theory}}`\Eqsheav{I_{\theory}}{M_{\theory}};\cat{M}_t`\cat{M}`v^*]
\dtriangle(0,300)/<-`<-`>/<750,300>[\Sets/M_{\theory}`\synt{C}{T}`\Sh{M_{\theory}};\cat{M}_d`u^*`]
%
%
%
%
\btriangle(1750,300)|ara|/->>`->>`->>/<600,300>[\Sets/M_{\theory}`\Sh{M_{\theory}}`\Sh{\synt{C}{T}};u`m_d`m_t]
\ptriangle(1750,0)/->>`->>`<<-/<600,300>[\Sh{M_{\theory}}`\Sh{\synt{C}{T}}`\Eqsheav{I_{\theory}}{M_{\theory}};`v`m]
\efig\]
where $u^*$ and $v^*$ are forgetful functors. The diagram on the right then shows the induced geometric morphisms. Showing that $\cat{M}$ is full and conservative, and that \synt{C}{T} generates \Eqsheav{I_{\theory}}{M_{\theory}} (as a full subcategory),
we conclude that $m$ is an equivalence:
\[\Eqsheav{I_{\theory}}{M_{\theory}}\ \simeq\classtop \]
\subsubsection{Sheaves on the space of models}
\label{Subsubsection: Sheaves on the set of models}
Given an object \syntob{\alg{x}}{\phi} of \synt{C}{T} we define the set
\[\sox{\alg{x}}{\phi}:=\cterm{\pair{\alg{M},\alg{[a]}}}{\alg{M}\in M_{\theory},\ \alg{[a]}\in \csem{\alg{x}}{\phi}^{\alg{M}}}\to^{\pi_{1}} M_{\theory}\]
over $M_{\theory}$, where $\pi_1$ projects out the model. Note the notation ``\sox{\alg{x}}{\phi}'' for the set on the left, which
we shall make extensive use of below. The mapping $\syntob{\alg{x}}{\phi}\mapsto
(\pi_1:\sox{\alg{x}}{\phi}\rightarrow M_{\theory})$ gives us the object part of a functor,
\[
\cat{M}_d:\synt{C}{T}\to \Sets/M_{\theory}
\]
(which sends an arrow of \synt{C}{T} to the obvious function over $M_{\theory}$).
\begin{Lemma}
\label{Lemma: Stone representation for coherent theories}
The functor
\[\cat{M}_d:\synt{C}{T}\to \Sets/M_{\theory}\]
is geometric and conservative, that is, $\cat{M}_d$ is faithful and reflects isomorphisms.
\begin{Proof}
Considering each \theory-model \alg{M} as a geometric functor from \synt{C}{T} to \Sets, we have a commuting triangle:
\[\bfig
\Atriangle/>`>`/<400,400>[\synt{C}{T}`\Sets/{M_{\theory}}`\prod_{\alg{M}\in
M_{\theory}}\Sets_{\alg{M}};\cat{M}_d`\pair{\ldots, \alg{M},\ldots}`]
\place(350,0)[\simeq]
\efig\]
Then $\cat{M}_d$ is geometric since all $\alg{M}\in M_{\theory}$ are geometric, and $\cat{M}_d$ is conservative
since the $\alg{M}\in M_{\theory}$ are enough for \theory\ (as models), and thus are jointly conservative (as functors).
\end{Proof}
\end{Lemma}
Next, we factor $\cat{M}_d$ through \Sh{M_{\theory}} by equipping the sets $\sox{\alg{x}}{\phi}$ with a topology such that the projection $\pi_1:\sox{\alg{x}}{\phi}\rightarrow M_{\theory}$ becomes a local homeomorphism.
\begin{Definition}\label{Definition: Logical topology on sheaves}
For an object $\syntob{\alg{x}}{\phi}$ of \synt{C}{T}, the \emph{logical topology} on the set
\[ \sox{\alg{x}}{\phi}=\cterm{\pair{\alg{M},\alg{[a]}}}{\alg{M}\in M_{\theory},
\alg{[a]}\in \csem{\alg{x}}{\phi}^{\alg{M}}} \] is the coarsest topology such that $\pi_1:\sox{\alg{x}}{\phi}\rightarrow M_{\theory}$ is continuous and such that for all lists $\alg{a}\in \thry{S}$ of the same length as \alg{x}. the image of the map
\[s_{\alg{a}}:\bopen{\syntob{\alg{x}}{\phi},\alg{a}}\to \sox{\alg{x}}{\phi}\]
defined by $\alg{M}\mapsto \pair{\alg{M},\alg{[a]}}$ is open.
\end{Definition}
There is an alternative characterization of this topology:
\begin{Lemma}\label{Lemma: Alt def of log top on sheaves}
Sets of the form
\[
\bopen{\syntob{\alg{x},\alg{y}}{\psi}, \alg{b}}:= \cterm{\pair{\alg{M},\alg{a}}}{ [\alg{a}] , [\alg{b}]\in
\csem{\alg{x},\alg{y}}{\phi\wedge\psi}^{\alg{M}}}
\]
where $\alg{b}$ is a tuple of elements of \thry{S} of the same length as $\alg{y}$, form a basis for the logical topology on \sox{\alg{x}}{\phi}.
\begin{Proof}
Straightforward.
\end{Proof}
\end{Lemma}
\begin{Remark}
Similar to Lemma \ref{Lemma: The logical topologies are all the same} specifying the topology on \sox{\alg{x}}{\phi} in terms of geometric formulas is redundant; Horn formulas will do.
\end{Remark}
For any object \syntob{\alg{x}}{\phi} in \synt{C}{T}, we now have the following:
\begin{Lemma}\label{Lemma: p is LH}
The projection $\pi_1:\sox{\alg{x}}{\phi}\rightarrow M_{\theory}$ is a local homeomorphism.
\begin{Proof}
It suffices to show that the projection is open. Using Lemma \ref{Lemma: Alt def of log top on sheaves}, let  $\bopen{\syntob{\alg{x},\alg{y}}{\psi}, \alg{b}} \subseteq \sox{\alg{x}}{\phi}$ be a basic open set. Then
\[
\pi_1\left(\bopen{\syntob{\alg{x},\alg{y}}{\psi}, \alg{b}}\right)= \bopen{\syntob{\alg{y}}{
\fins{\alg{x}}\phi \wedge \psi}, \alg{b}}\subseteq M_{\theory}
\]
which is open.
\end{Proof}
\end{Lemma}
\begin{Lemma}\label{Lemma: maps are continuous} Given an arrow
\[\sigma:\syntob{\alg{x}}{\phi}\to
\syntob{\alg{y}}{\psi}\]
in \synt{C}{T}, the corresponding function $\cat{M}_d(\sigma): \sox{\alg{x}}{\phi}\rightarrow \sox{\alg{y}}{\psi}$ is
continuous.
\begin{Proof} Given a basic open
$\bopen{\syntob{\alg{y},\alg{z}}{\xi}, \alg{c}} \subseteq \sox{\alg{y}}{\psi}$, then
\[
\cat{M}_d(\sigma)^{-1}\left(\bopen{\syntob{\alg{y},\alg{z}}{\xi}, \alg{c}} \right) =
\bopen{\syntob{\alg{x},\alg{z}}{\fins{\alg{y}} \sigma \wedge \xi}, \alg{c}}
\]
\end{Proof}
\end{Lemma}
%
%
%
We conclude with the following proposition (see also Theorem 6.3.3 of \cite{makkaireyes}).
\begin{Proposition}\label{Proposition: DC, Md factors as M}
The functor
$\cat{M}_d:\synt{C}{T}\to\Sets/M_{\theory}$
factors through the category \Sh{M_{\theory}} of sheaves as
\[
\bfig \dtriangle/<-`<-`->/<600,400>[\Sets/M_{\theory}`\synt{C}{T}`\Sh{M_{\theory}};\cat{M}_d`u^*`\cat{M}_t] \efig
\]
where $u^*:\Sh{M_{\theory}}\to\Sets/M_{\theory}$ is the forgetful (inverse image) functor. Moreover, $\cat{M}_t$ is
geometric and conservative.
\begin{Proof}$\cat{M}_t$ is obtained by Lemma \ref{Lemma: p is LH} and Lemma \ref{Lemma: maps are
continuous}.
Since $\cat{M}_d$ is geometric and the forgetful functor $u^*$ reflects geometric structure (being geometric and conservative), $\cat{M}_t$ is
geometric. Since $\cat{M}_d$ is conservative, $\cat{M}_t$ is, too.
%
%
%
\end{Proof}
\end{Proposition}
\subsubsection{Equivariant sheaves on the space of models}
\label{Subsubsection: Equivariant sheaves on the space of models}
Consider the groupoid of \theory-models and isomorphisms  $\topo{M}_{\theory}=I_{\theory}\rightrightarrows M_{\theory}$. For an object \syntob{\alg{x}}{\phi} in
\synt{C}{T}, we have the functor $\cat{M}_t:\synt{C}{T}\to \Sh{M_{\theory}}$ assigning \syntob{\alg{x}}{\phi} to the sheaf $\pi_1:\sox{\alg{x}}{\phi}\rightarrow M_{\theory}$. There is an obvious action of `application',
\begin{equation}\label{Equation: The application action}
\theta_{\syntob{\alg{x}}{\phi}}:I_{\theory}\times_{M_{\theory}} \sox{\alg{x}}{\phi}\to \sox{\alg{x}}{\phi}
\end{equation}
defined by sending a \theory-isomorphism $\alg{f}:\alg{M}\to \alg{N}$ and an element $\pair{\alg{M},\alg{[a]}}\in
\csem{\alg{x}}{\phi}^{\alg{M}}$ to the element $\pair{\alg{N},\alg{f}(\alg{[a]})}\in \csem{\alg{x}}{\phi}^{\alg{N}}$. We shall mostly
leave the subscript implicit. When it is clear from context in which fiber over $M_{\theory}$ an element of \sox{\alg{x}}{\phi} lies, we shall often also leave this implicit and write, say, $[\alg{a}]$ instead of $\pair{\alg{M},[\alg{a}]}$.

\begin{Lemma}\label{Lemma: Definable set functors are equivariant sheaves}
The sheaf $\pi_1:\sox{\alg{x}}{\phi}\rightarrow M_{\theory}$ together with the function $\theta:I_{\theory}\times_{M_{\theory}}
\sox{\alg{x}}{\phi}\to \sox{\alg{x}}{\phi}$ is an object of \Eqsheav{I_{\theory}}{{M_{\theory}}}.
\begin{Proof}
We must verify that $\theta$ is continuous and satisfies the axioms for being an action. The latter is
straightforward, so we do the former.  Let $V=\bopen{\syntob{\alg{x},\alg{y}}{\psi},\alg{a}}\subseteq
\sox{\alg{x}}{\phi}$ be given, and suppose $\pair{f:\alg{M}\rightarrow \alg{N},\alg{[b]}}\in
\theta^{-1}(V)\subseteq I_{\theory}\times_{M_{\theory}} \sox{\alg{x}}{\phi}$. Choose $\alg{c}$ such that
$f(\alg{[c]})=\alg{[a]}$. Then
\[
\pair{f:\alg{M}\rightarrow \alg{N},\alg{[b]}} \in W:= \left(\begin{array}{c}
-  \\
\alg{c} \mapsto \alg{a}   \\
 -
\end{array}\right)  \times_{M_{\theory}}
\bopen{\syntob{\alg{x},\alg{y}}{\psi},\alg{c}} \]
and $\theta(W)\subseteq V$. So $\theta$ is continuous.
\end{Proof}
\end{Lemma}
It is clear that any definable morphism of sheaves
\[\cat{M}_t(\sigma):\cat{M}_t(\syntob{\alg{x}}{\phi})\to \cat{M}_t(\syntob{\alg{y}}{\psi})\]
for $\sigma: \syntob{\alg{x}}{\phi} \to \syntob{\alg{y}}{\psi}$ in \synt{C}{T},
commutes with the respective actions of application, and so we have a functor
\[ \cat{M}:\synt{C}{T}\to \Eqsheav{I_{\theory}}{{M_{\theory}}} \]
The faithful forgetful functor $v^*:\Eqsheav{I_{\theory}}{{M_{\theory}}}\to \Sh{{M_{\theory}}}$ is the inverse image part of a
geometric morphism $v:\Sh{{M_{\theory}}}\epi \Eqsheav{I_{\theory}}{{M_{\theory}}}$. Composing $\cat{M}$ with $v^*$ we get a commuting triangle:
\[ \bfig \dtriangle/<-`<-`>/<450,300>[\Sh{{M_{\theory}}}`\synt{C}{T}`\Eqsheav{I_{\theory}}{{M_{\theory}}};\cat{M}_t`v^*`\cat{M}]   \efig \]
from which we conclude that $\cat{M}$ is geometric and conservative, and that we have a factorization,
\[ \bfig \dtriangle/->>`->>`<<-/<550,300>[\Sh{{M_{\theory}}}`\classtop`\Eqsheav{I_{\theory}}{{M_{\theory}}};m_t`v`m]   \efig \]
We state these facts for reference:
\begin{Lemma}\label{Lemma: m dager surjective}
$\cat{M}:\synt{C}{T}\to\Eqsheav{I_{\theory}}{M_{\theory}}$ is geometric and conservative (i.e.\ faithful and
reflects isomorphisms).
\end{Lemma}
Next, we aim to show that the geometric morphism $m:\Eqsheav{I_{\theory}}{{M_{\theory}}}\epi\classtop$ is an equivalence by showing that
$\cat{M}(\synt{C}{T})$ is a site for \Eqsheav{I_{\theory}}{{M_{\theory}}}. First, it is clear that subobjects of an equivariant sheaf
\pair{a:A\rightarrow {M_{\theory}}, \alpha} can be thought of, and represented, as open subsets of $A$ that are closed under
the action $\alpha$. We call a subset $S\subseteq A$ that is
closed under the action \emph{stable} (in order to reserve \emph{closed} to mean topologically closed), and we call the least stable subset containing $S$ the
\emph{stabilization} of $A$. We call the objects and arrows in the image of $\cat{M}:\synt{C}{T}\to
\Eqsheav{I_{\theory}}{{M_{\theory}}}$ \emph{definable}.
\begin{Lemma}\label{Lemma: Basic opens stabilize as definables}
Let \alg{a} be a tuple of distinct elements of \thry{S}. The stabilization of a basic open subset \[\bopen{\syntob{\alg{x},\alg{y}}{\psi},\alg{a}} \subseteq
\sox{\alg{x}}{\phi}\] is the definable subset $\sox{\alg{x}}{\phi\wedge \fins{\alg{y}}\psi}\subseteq
\sox{\alg{x}}{\phi}$.
\begin{Proof}
Let $\pair{\alg{M},\alg{[b]}}\in \sox{\alg{x}}{\phi\wedge \fins{\alg{y}}\psi}$ be given, and choose a witness \alg{[c]} in \alg{M} for the existential quantifier. By Lemma \ref{Lemma: Star lemma} there is a model \alg{N} and an isomorphism $\alg{f}:\alg{N}\rightarrow\alg{M}$ such that $\alg{f}(\alg{[a]})=\alg{[c]}$. Then $\alg{f}^{-1}(\alg{[b]}),\alg{[a]}\in \csem{\alg{x},\alg{y}}{\phi\wedge\psi}^{\alg{N}}$, so that $\pair{\alg{N},\alg{f}^{-1}(\alg{[b]})}\in \bopen{\syntob{\alg{x},\alg{y}}{\psi},\alg{a}} \subseteq\sox{\alg{x}}{\phi}$ and $\theta(\alg{f}:\alg{N}\rightarrow \alg{M},\alg{f}^{-1}(\alg{[b]}))=\pair{\alg{M},\alg{[b]}}$.
\end{Proof}
\end{Lemma}
\begin{Corollary}\label{Corollary: Mdag full on subobjects}
The functor $\cat{M}:\synt{C}{T}\to \Eqsheav{I_{\theory}}{{M_{\theory}}}$ is full on subobjects.
\begin{Proof} Stabilizing commutes with unions.
\end{Proof}
\end{Corollary}
\begin{Corollary}\label{Corollary: Mdag is full}
The functor $\cat{M}:\synt{C}{T}\to \Eqsheav{I_{\theory}}{{M_{\theory}}}$ is full.
\begin{Proof} Since $\cat{M}$ is full on subobjects, geometric, and conservative, any functional relation between two objects in $\cat{M}(\synt{C}{T})$ comes from a functional relation in \synt{C}{T}.
\end{Proof}
\end{Corollary}
\subsubsection{Semantic representation of \theory}
\label{Subsubsection: Semantic representation of T}
The following sequence of lemmas serve to establish that the definable objects form a generating set for \Eqsheav{I_{\theory}}{M_{\theory}}, thus combining with Corollary \ref{Corollary: Mdag is full} to show that $\Eqsheav{I_{\theory}}{M_{\theory}}\simeq \classtop$. We show that the definables are a dense subcategory (in the sense of \cite[C2.2.1]{elephant1}) of the Moerdijk site described in \ref{Subsection: Equivariant sheaves on topological groupoids} (a direct proof can also be given, but it is somewhat tedious and, perhaps, of less interest). Note (again) that given a basic open $\bopen{\syntob{\alg{x}}{\phi},\alg{a}}\subseteq M_{\theory}$, we can always assume without loss of generality that \alg{a} is a tuple of distinct elements of \thry{S} (or we can shorten the context and rewrite $\phi$ accordingly).

Given \syntob{\alg{x}}{\phi} and a list \alg{a} of distinct elements of \thry{S}, the section $s_{\alg{a}}:\bopen{\syntob{\alg{x}}{\phi},\alg{a}}\rightarrow \sox{\alg{x}}{\phi}$ defined by
\[\alg{M}\mapsto \pair{\alg{M},\alg{[a]}}\in \sox{\alg{x}}{\phi}\]
is continuous by Definition \ref{Definition: Logical topology on sheaves} and Lemma \ref{Lemma: p is LH}.
\begin{Lemma}\label{Lemma: Basic opens of X have stabilizers}
Let a basic open $\bopen{\syntob{\alg{x}}{\phi},\alg{a}}\subseteq M_{\theory}$ be given, where we assume without loss that the $a$'s are
distinct. Then $\sox{\alg{x}}{\phi}$ is the stabilization of
$s_{\alg{a}}(\bopen{\syntob{\alg{x}}{\phi},\alg{a}})$ (with respect to the action of application $\theta$).
\begin{Proof}
By Lemma \ref{Lemma: Basic opens stabilize as
definables}, $\sox{\alg{x}}{\phi}$ is the stabilization of
$s_{\alg{a}}(\bopen{\syntob{\alg{x}}{\phi},\alg{a}})=\bopen{\syntob{\alg{x},\alg{y}}{\alg{x}=\alg{y}},\alg{a}}$.
\end{Proof}
\end{Lemma}
%
%
\begin{Lemma}\label{Lemma: Definables are GUNS}
Let $\syntob{\alg{x}}{\phi}\in\synt{C}{T}$ and let \alg{a} be a list of distinct elements of \thry{S} of the same length as \alg{x}. Then the section $s:\bopen{\syntob{\alg{x}}{\phi},\alg{a}}\rightarrow
\sox{\alg{x}}{\phi}$ defined by $\alg{M}\mapsto \pair{\alg{M},\alg{[a]}}$  induces an isomorphism
\[\hat{s}:\pair{\topo{M}_{\theory},\bopen{\syntob{\alg{x}}{\phi},\alg{a}},d^{-1}(\bopen{\syntob{\alg{x}}{\phi},\alg{a}})\cap\bopen{\alg{a}\mapsto\alg{a}}}\to \cat{M}(\syntob{\alg{x}}{\phi})\]
\begin{Proof}
Let $U=\bopen{\syntob{\alg{x}}{\phi},\alg{a}}$. We have
\[N_s=\cterm{\alg{f}\in d^{-1}(U)}{\theta(\alg{f},s(d(\alg{f})))=s(c(\alg{f}))}=d^{-1}(\bopen{\syntob{\alg{x}}{\phi},\alg{a}})\cap\bopen{\alg{a}\mapsto\alg{a}}\]
whence $s$ lifts to a morphism $\hat{s}:\pair{\topo{M}_{\theory},U,N_s}\rightarrow \cat{M}(\syntob{\alg{x}}{\phi})$ defined by $\hat{s}(\alg{g})=\theta(\alg{g},s(d(\alg{g})))$. This morphism is easily verified to be injective, and it is surjective by Lemma \ref{Lemma: Basic opens of X have stabilizers}.
\end{Proof}
\end{Lemma}
\begin{Lemma}\label{Lemma: Finding definable in GUN}
Let an open subset $N\subseteq I_{\theory}$ closed under inverses and composition be given, and let $\alg{M}\in d(N)=c(N)=U$. Then there exist a $\syntob{\alg{x}}{\phi}\in\synt{C}{T}$ and a list of distinct elements \alg{a} of \thry{S} of the same length as \alg{x} such that
\[1_{\alg{M}}\in \left(\begin{array}{c}
\syntob{\alg{x}}{\phi},\alg{a}  \\
\alg{a} \mapsto \alg{a}   \\
 \syntob{\alg{x}}{\phi},\alg{a}
\end{array}\right)\subseteq N\]
\begin{Proof}
 Choose a basic open neighborhood $V$ of $1_{\alg{M}}$ in $N$. Informally explained, $V$ can just be `rewritten' to obtain a smaller open neighborhood of the required form which still contains $1_{\alg{M}}$: domain and codomain conditions can be added to each other; elements of \thry{S} occurring in the domain or codomain condition can be added to the preservation condition; and if the preservation condition contains $b\mapsto c$ for distinct $b$ and $c$, then it can be replaced by $b\mapsto b$ and the domain and codomain condition $b=c$.
\end{Proof}
\end{Lemma}
\begin{Proposition}\label{Lemma: The definables generate equivariant sheaves on the groupoid of models}
The category of definable objects is a dense subcategory of the Moerdijk site of \Eqsheav{I_{\theory}}{M_{\theory}}. In particular, the set of definable objects in \Eqsheav{I_{\theory}}{M_{\theory}} is  generating.
\begin{Proof}Let \pair{\topo{M}_{\theory},U,N} be given, and let $[\alg{f}]\in d^{-1}(U)/_{\sim_N}$. It suffices to show that $[\alg{f}]$ is in the image of a morphism with definable domain.  Say $\alg{M}=d(\alg{\alg{f}})$. By Lemma \ref{Lemma: Finding definable in GUN}, we can find
\[1_{\alg{M}}\in M=\left(\begin{array}{c}
\syntob{\alg{x}}{\phi},\alg{a}  \\
\alg{a} \mapsto \alg{a}   \\
 \syntob{\alg{x}}{\phi},\alg{a}
\end{array}\right)\subseteq N\]
Let $V=\bopen{\syntob{\alg{x}}{\phi},\alg{a}}$. The continuous section $e:V\rightarrow d^{-1}(U)/_{\sim_N}$ defined by $\alg{K}\mapsto [1_{\alg{K}}]$ lifts to a morphism $\hat{e}:\pair{\topo{M}_{\theory},V,M}\rightarrow \pair{\topo{M}_{\theory},U,N}$ such that $\hat{e}[g]_{\sim_M}=[g]_{\sim_N}$ (see \cite{moerdijk:88} or the expose in \cite{forssell:subgpds}). And so $[\alg{f}]$ is in the image of $\hat{e}$. And by Lemma \ref{Lemma: Definables are GUNS}, $\pair{\topo{M}_{\theory},V,M}\cong\cat{M}(\syntob{\alg{x}}{\phi})$.
\end{Proof}
\end{Proposition}
We can now conclude with the following representation result, which, minor differences aside, should be attributed to Butz and Moerdijk (see \cite{butz:98b}).
\begin{Theorem}\label{Theorem: The classifying topos of T is equivariant G-sheaves}
For \theory\ a geometric theory with enough \thry{S}-indexed models there is an equivalence
\[ \classtop\simeq \Eqsheav{I_{\theory}}{{M_{\theory}}} \]
where $\topo{M}_{\theory} = \left(I_{\theory}\rightrightarrows M_{\theory}\right)$ is the topological groupoid
of \thry{S}-indexed \theory-models.
\begin{Proof}
By Proposition \ref{Lemma: The definables generate equivariant sheaves on the groupoid of models}, the definable
objects form a generating set for \Eqsheav{I_{\theory}}{{M_{\theory}}}. Therefore, the full subcategory of definable objects
equipped with the coverage inherited from the canonical coverage of \Eqsheav{I_{\theory}}{{M_{\theory}}} is a site for
\Eqsheav{I_{\theory}}{{M_{\theory}}} (see \cite[C2.2.16]{elephant1}). By Lemma \ref{Lemma: m dager surjective} and Corollary \ref{Corollary: Mdag is full}, the full subcategory of definable
objects is equivalent to
\synt{C}{T}, and the canonical coverage inherited from \Eqsheav{I_{\theory}}{{M_{\theory}}} is just the geometric coverage, $G$, on \synt{C}{T}. Therefore, \pair{\synt{C}{T}, G} is a site for
\Eqsheav{I_{\theory}}{{M_{\theory}}}. But the topos of sheaves on \pair{\synt{C}{T}, G} is the classifying topos of \theory\ (see \cite{elephant1}),  whence \[\Eqsheav{I_{\theory}}{{M_{\theory}}}\simeq \Sh{\synt{C}{T},G}\simeq \classtop\]
\end{Proof}
\end{Theorem}
\section{Syntax-semantics duality for geometric theories}
\label{Section: Syntax-semantics duality}
Theorem \ref{Theorem: The classifying topos of T is equivariant G-sheaves} tells us that the classifying topos of a geometric theory with enough models can be constructed both syntactically (from the syntactic category) and semantically from the groupoid $\topo{M}_{\theory}$ of models and isomorphisms, constructed using a sufficiently large `index' set \thry{S}. Accordingly, toposes with enough points can be regarded from both a logical, syntactical point of view and from a geometric semantic point of view. We use this to give a duality between the category of (syntactical categories of single-sorted) geometric theories and a slice category of topological groupoids, in the form of a `syntax-semantics' adjunction the counit components of which at sufficiently small theories with enough models are isomorphisms. Fix an index set \thry{S} of size $\kappa$, for some infinite cardinal. The adjunction is constructed using \thry{S} to construct semantical groupoids, and, accordingly, theories are `sufficiently small' if they have enough models of size smaller or equal to $\kappa$. Conversely, any topological groupoid `over \thry{S}' (in a sense to be specified in Section \ref{Subsection: Singlesort, The Object Classifier}) gives rise to such a theory.

\subsection{The category of theories}
\label{Subsection: Theories and Categories}
In order to recover a theory from a groupoid `over \thry{S}' in a canonical way, we make the following specifications.
Let \theory\ be a geometric theory over a first-order single-sorted signature. Recall that we construct the \emph{syntactic category} \synt{C}{T} of a single-sorted theory \theory\ as
follows. The objects of \synt{C}{T} are equivalence classes of ($\alpha$-equivalence classes of)
formulas-in-context, $|\syntob{\alg{x}}{\phi}|$, where $\syntob{\alg{x}}{\phi}\sim \syntob{\alg{x}}{\psi}$ iff
$\theory$ proves the sequents $ \phi\dashv\vdash_{\alg{x}} \psi$. Arrows between such objects are as usual given by
\theory-provable equivalence classes of formulas-in-context,
\[ |\syntob{\alg{x},\alg{y}}{\sigma}|:|\syntob{\alg{x}}{\phi}|\to |\syntob{\alg{y}}{\psi}|\]
such that $\theory$ proves that $\sigma$ is a functional relation between $\phi$ and $\psi$. This definition  of \synt{C}{T} is clearly equivalent, in the sense of
producing equivalent categories, to the usual one where objects are just $\alpha$-equivalence classes (and not
\theory-provable equivalence classes) of formulas-in-context, but is more convenient as long as we are mostly
interested in \theory-models in \Sets. Moreover, it results in a small syntactic category (see \cite[D1.3.8]{elephant1}). In what follows, we usually drop the vertical bars indicating equivalence
class in our notation (i.e.\ we write \syntob{\alg{x}}{\phi} but mean $|\syntob{\alg{x}}{\phi}|$). With this
definition of syntactical category, every syntactic category has the properties:
\begin{itemize}
\item There is a  distinguished object, $U$,  with distinguished distinct finite powers.
(In a syntactic category $U=\syntob{x}{\top}$.)
\item There is a \emph{system of inclusions}, that is, a
set $\mathfrak{I}$ of distinguished monomorphisms which is closed under composition
and  identities, and such that every object has a unique inclusion into a finite power of $U$. Moreover (and
this is not the case with the alternative definition of \synt{C}{T})
 every subobject, considered as a set of monomorphisms, of an object contains a unique inclusion. (We can take
 $\mathfrak{I}$ in \synt{C}{T} to be
  the set of all arrows $\syntob{x_1,\ldots,x_n}{\phi}\to\syntob{y_1,\ldots,y_n}{\psi}$ which contain the
 formula-in-context \syntob{\alg{x},\alg{y}}{\phi \wedge \psi \wedge \alg{x}=\alg{y}}.)
\end{itemize}
We claim that this characterizes syntactical categories (for single-sorted  theories) up to
isomorphism. Suppose \cat{B} is a (small) geometric category (see \cite[A1.4.18, D1.4]{elephant1}) with a  distinguished object and a system of
inclusions. Let the signature, $\Sigma_{B}$, of \cat{B} consist of,  for each inclusion $R\embedd U^n$ a n-ary
relation symbol. Set the theory,  $\theory_{\cat{B}}$, of \cat{B} to be the set of true geometric sequences over $\Sigma_{B}$
under the canonical interpretation in \cat{B}.
\begin{Lemma}\label{Lemma: Characterizing syntcats}
There is an isomorphism
\[\cat{B}\cong \synt{C}{\theory_{\cat{B}}}.\]
\begin{Proof}
Define a functor $F:\cat{B}\to\synt{C}{T}$ by sending an object $A$ in \cat{B} to
\syntob{\alg{x}}{R_A}, where $R_A$ is the predicate in $\Sigma_{B}$ corresponding to its unique inclusion
$A\embedd U^n$ into a power of $U$. For an arrow $f:A\to B$ in \cat{B},  there is
an inclusion $\mng{Grph}(f)\embedd U^{n+m}$ corresponding to a relation symbol $R_{\mng{Grph}(f)}$ such that
\syntob{\alg{x},\alg{y}}{R_{\mng{Grph}(f)}} is $\theory_{\cat{B}}$-provably functional from \syntob{\alg{x}}{R_A} to
\syntob{\alg{y}}{R_B}, so set $F(f)=\syntob{\alg{x}, \alg{y}}{R_{\mng{Grph}(f)}}$. In the other direction, define a
functor $G:\synt{C}{T}\to\cat{B}$ by sending an object \syntob{\alg{x}}{\phi} to the domain of the inclusion
representing the subobject $\csem{\alg{x}}{\phi}^{\cat{B}}$ under the canonical interpretation of $\Sigma_{B}$ in
\cat{B}. Then $G\circ F = 1_{\cat{B}}$. And if $\syntob{\alg{x}}{\phi} \in \synt{C}{T}$, with
$F(\syntob{\alg{x}}{\phi})=A\embedd U^n$, and \syntob{\alg{x}}{R_A} is the predicate of $\Sigma_{B}$ representing
$A\embedd U^n$, then $\theory_{\cat{B}}$ proves $ \phi \dashv \vdash_{\alg{x}}  R_A$, so $F\circ G=
1_{\synt{C}{T}}$.
\end{Proof}
\end{Lemma}
Note that if $F:\cat{B} \to \cat{D}$ is a
geometric functor that preserves the  distinguished object, or synonymously the \emph{single sort}, then it is
naturally isomorphic to one that moreover preserves the distinguished finite powers of $U$ on the nose, and that
preserves inclusions.
\begin{Definition}\label{Definition: Singlesort, FOL}
The category \cat{T} consists of geometric categories with a distinguished object and a system of inclusions. Arrows in \cat{T} are geometric
functors that preserve the distinguished object (and its distinguished finite powers) and inclusions on the
nose.
\end{Definition}
We write \synt{C}{T} for an object of \cat{T}, since it is (isomorphic to) a syntactic category for a
geometric theory \theory\ by Lemma \ref{Lemma: Characterizing syntcats}. By a \thry{T}-model, we mean a geometric
functor $M:\synt{C}{T}\to\Sets$ that sends $M(U^n)$ to the n-fold cartesian product of $M(U)$, and inclusions to
subset inclusions.
\subsection{The object classifier}
\label{Subsection: Singlesort, The Object Classifier}
%
%
%
%
Denote by $\theory_{=}$ the (single-sorted) geometric theory with no constant, function, or relation symbols (except equality) and no axioms. Accordingly, $\mathbf{Set}[\theory_{=}]$ classifies objects in the category
of toposes and geometric morphisms. Since in the
category \cat{T} there are only distinguished object preserving functors, there exists exactly one
arrow from \synt{C}{T_=} to any \synt{C}{T} in \cat{T}, that is, \synt{C}{T_=} is an initial object. Dually, we consider topological groupoids over the semantical groupoid $\topo{M}_{\theory_=}$.
\begin{equation}\label{Equation: Functor triangles}\bfig
\Vtriangle/>`<-`<-/<350,350>[\synt{C}{T}`\synt{C}{T'}`\synt{C}{T_=};``]
\place(1000,175)[\mapsto]
\Vtriangle(1500,0)/<-`>`>/<350,350>[\topo{M}_{\thry{T}}`\topo{M}_{\thry{T'}}`\topo{M}_{\thry{T_=}};``]
\efig\end{equation}
where $\topo{M}_{\theory}$ is the topological groupoid of models (with underlying set a quotient of \thry{S}) and isomorphisms equipped with the logical topology (Definition \ref{Definition> Logical topology}). Since \thry{S} is at least countable and $\theory_=$ has enough countable models, Theorem \ref{Theorem: The classifying topos of T is equivariant G-sheaves} applies and $\mathbf{Set}[\theory_{=}]\simeq \Eqsheav{I_{\theory_=}}{M_{\theory_=}}$. Given the importance of $\topo{M}_{\theory_=}$ in the construction of the syntax-semantics adjunction below, we introduce new notation for it and spell out what it consists of:
\begin{Definition}\label{Definition: Singlesort, N}
The topological groupoid $\topo{S}$ consists of the set $S_0$ of quotients of subsets of \thry{S} with the set
$S_1$ of bijections between them, equipped with topology as follows. The topology on the set of objects
%
%
%
is the coarsest topology in which sets of the form
\[
(a\sim b):=\cterm{A\in S_0}{[a]=[b]\ \textnormal{in}\ A}
\]
are open. As in Definition \ref{Definition> Logical topology}, this condition should be taken to mean that $A$ is a quotient of a subset which contains $a$ and $b$ and that $a$ and $b$ are members of the same equivalence class in $A$.   The topology on the set, $S_1$ of bijections is the coarsest topology such
that the domain and codomain maps $d,c:S_1\rightrightarrows S_0$ are both continuous, and such that all sets of
the form
\[
(a\mapsto b):=\cterm{f:A\to<150>^{\cong} B\ \textnormal{in}\ S_1}{ f([a])=[b]}
\]
are open.
\end{Definition}
Comparing with Definition \ref{Definition> Logical topology}, we see that this simply restates the logical topology for the signature which only contains equality, so that
\[\topo{S}\cong  \topo{M}_{\theory_=}\]
The category \Eqsheav{S_1}{S_0} of equivariant sheaves on $\topo{S}$, therefore, classifies
objects. The generic  object, $\cat{U}$, in \Eqsheav{S_1}{S_0} can be
taken to be the definable sheaf
\[\pair{\sox{x}{\top}\rightarrow M_{\theory_{=}}\cong
S_0,\theta_{\syntob{x}{\top}}}\]
which we can characterize directly by restating Definition \ref{Definition: Logical topology on sheaves}:
\begin{Lemma}\label{Definition: U generic supported object}
The generic  object, $\cat{M}(\syntob{x}{\top})$ in \Eqsheav{S_1}{S_0} can be characterized as the equivariant sheaf
\[\cat{U}=\pair{p:U\rightarrow S_0, \theta}\]
where $U=\coprod_{A\in S_0}A$; the function $p:U\rightarrow S_0$ is the projection; the topology on $U$ is the coarsest such that $p$ is continuous and every set of the form \[\bopen{a}:=\cterm{\pair{A,[a]}}{A\in S_0\ \textnormal{is a quotient of a subset of \thry{S} containing $a$}}\] is open; and $\theta$ is the obvious action.
\end{Lemma}
\begin{Remark}Since \Eqsheav{S_1}{S_0} classifies objects, with \cat{U} being the generic object, it is equivalent (see e.g.\ \cite{maclane92}) to the functor category of sets to finite sets
\[\Eqsheav{S_1}{S_0}\simeq \Sets^{\alg{Fin}}\]
in which the generic object can be taken to be the inclusion $\alg{Fin}\hookrightarrow\Sets$.
\end{Remark}
Since \Eqsheav{S_1}{S_0} classifies objects, a morphism $f:\topo{G}\to\topo{S}$ of topological groupoids induces a geometric morphism $f:\Eqsheav{G_1}{G_0}\to\Eqsheav{S_1}{S_0}$ which classifies an  object $f^*(\cat{U})$ in \Eqsheav{G_1}{G_0}. This allows us to construct an adjunction between \cat{T}
and the category $\alg{Gpd}/\topo{S}$ of topological groupoids over \topo{S}.
\begin{Remark}\label{Remark: Reindexing} Throughout this paper we will stick with a fixed index set \thry{S} and leave out the issue of switching index sets. In passing we briefly note, however, that switching to  larger index sets (larger initial segments of the set theoretic universe, say) induces the expected morphisms. If $\thry{S}\subseteq \overline{\thry{S}}$, the topological groupoid $\overline{\topo{S}}$ can be defined from $\overline{\thry{S}}$ in the same way that ${\topo{S}}$ is defined from ${\thry{S}}$, and the inclusion induces a morphism of topological groupoids,
\[i:\topo{S}\to\overline{\topo{S}}\]
with subspace inclusions as component functions (one easily sees for instance that $S_0$ is a closed subspace of $\overline{S_0}$). Composing along $i$ yields a functor $i_*:\alg{Gpd}/\topo{S}\to \alg{Gpd}/\overline{\topo{S}}$, and restricting along $i$ is readily seen to be a right adjoint
\[\bfig
\morphism|b|/{@{<-}@/_1em/}/<1000,0>[\alg{Gpd}/\topo{S}`\alg{Gpd}/\overline{\topo{S}};i^*]
\morphism/{@{>}@/^1em/}/<1000,0>[\alg{Gpd}/\topo{S}`\alg{Gpd}/\overline{\topo{S}};i_*]
\place(500,0)[\bot]
\efig\]
\end{Remark}
\subsection{The semantical groupoid functor \Mod}
\label{Subsection: The Semantical Groupoid Functor}
For \theory\ a single-sorted geometric theory, regardless of whether it has enough \thry{S}-indexed models, we can construct the (possibly empty)
topological groupoid, $\topo{M}_{\theory}$, of \theory-models with underlying set a quotient of a subset of \thry{S} and
\theory-isomorphisms between them, equipped with the logical topology as in Section \ref{Section: Topological groupoids and categories of models}. This gives the object part of a functor $\cat{T}\to\alg{Gpd}$.
For a morphism $F:\synt{C}{T}\to \synt{C}{T'}$ in \cat{T}, any \thry{S}-model,
$\synt{C}{T'}\to\Sets$, `restricts' along $F$ to a \theory-model with the same underlying set, and any
\thry{T'}-model isomorphism restricts to a \theory-model isomorphism with the same underlying function, and so we get
functions $f_0:M_{\thry{T'}} \rightarrow M_{\theory}$ and $f_1:I_{\thry{T'}} \rightarrow I_{\theory}$ such that
the following commutes:
\[
\bfig
%
%
\morphism<750,0>[I_{\theory}\times_{M_{\theory}} I_{\theory}`I_{\theory};m]
%
%
\morphism(750,0)|a|/@{>}@<5pt>/<750,0>[I_{\theory}`M_{\theory};d]
\morphism(750,0)|m|/@{<-}/<750,0>[I_{\theory}`M_{\theory};e]
\morphism(750,0)|b|/@{>}@<-5pt>/<750,0>[I_{\theory}`M_{\theory};c]
\Loop(750,0)I_{\theory}(ur,ul)_{i}
%
%
%
\morphism(0,-500)<750,0>[I_{\thry{T'}}\times_{M_{\thry{T'}}} I_{\thry{T'}}`I_{\thry{T'}};m]
%
%
\morphism(750,-500)|a|/@{>}@<5pt>/<750,0>[I_{\thry{T'}}`M_{\thry{T'}};d]
\morphism(750,-500)|m|/@{<-}/<750,0>[I_{\thry{T'}}`M_{\thry{T'}};e]
\morphism(750,-500)|b|/@{>}@<-5pt>/<750,0>[I_{\thry{T'}}`M_{\thry{T'}};c]
\Loop(750,-500)I_{\thry{T'}}(dr,dl)_{i}
%
\morphism(0,-500)<0,500>[I_{\thry{T'}}\times_{M_{\thry{T'}}} I_{\thry{T'}}`I_{\theory}\times_{M_{\theory}}
I_{\theory};f_1\times f_1]
\morphism(750,-500)<0,500>[I_{\thry{T'}}`I_{\theory};f_1]
\morphism(1500,-500)<0,500>[M_{\thry{T'}}`M_{\theory};f_0]
\efig
\]
Given a basic open $\bopen{\syntob{\alg{x}}{\phi},\alg{a}}\subseteq M_{\theory}$, we see that its inverse image is given by
translating $\phi$ along $F$,
\[
f_0^{-1}\left( \bopen{\syntob{\alg{x}}{\phi},\alg{a}} \right) = \bopen{\syntob{\alg{x}}{F(\phi)},\alg{a}}\subseteq M_{\thry{T'}}
\]
(recall that $F$ preserves the finite powers of the distinguished object as well as inclusions so we can allow ourselves to write $F(\syntob{\alg{x}}{\phi})=\syntob{\alg{x}}{F(\phi)}$) and so $f_0$ is continuous. Similarly, $f_1$ is continuous because $f^{-1}(\bopen{a\mapsto b})=\bopen{a\mapsto b}$ Thus we obtain a morphism of continuous
groupoids $f:\topo{M}_{\thry{T'}}\to \topo{M}_{\theory}$, which is, then, the morphism-part of what is clearly a functor $\cat{T}\to\alg{Gpd}$.
To construct the `semantic' functor, we apply this functor to triangles of the form seen in Diagram (\ref{Equation: Functor triangles}) to obtain a functor from \cat{T} into the category of topological groupoids over \topo{S}. Specifically, a category \synt{C}{T} in \cat{T}, which has a unique morphism $U_{\theory}:\synt{C}{T_=}\to\synt{C}{T}$, is sent to the groupoid morphism $u_{\theory}:\topo{M}_{\theory}\to\topo{M}_{\theory_=}\cong\topo{S}$. The morphism $u_{\theory}$, then, is the forgetful morphism sending  models to their underlying sets.
\begin{Definition}\label{Definition: Gamma}
The contravariant \emph{semantical groupoid functor} \Mod
\[
\Mod:\cat{T}^{\text{op}} \to \alg{Gpd}/\topo{S}
\]
sends a theory to its topological groupoid over $\topo{S}$ of models:
\[
\bfig
\morphism(-500,-150)<400,0>[\synt{C}{T}`\synt{C}{T'};F]
\place(500,-150)[\mapsto]
\Vtriangle(1000,-300)<300,300>[\topo{M}_{\thry{T'}}`\topo{M}_{\thry{T}}`\topo{S};f`u_{\thry{T'}}`u_{\thry{T}}] \efig
\]
\end{Definition}
\begin{Remark}
As long as there is little danger of confusion,  we will continue the practice begun above of denoting the semantical
morphism of groupoids obtained by applying the semantical groupoid functor to a morphism of \cat{T} simply by
switching from an upper case to a lower case letter (and then using the same lower case letter for the induced geometric morphism). E.g.\ for $F:\synt{C}{T}\to\synt{C}{T'}$, we get the morphism of topological groupoids
\[ \Mod(F)=f:\topo{M}_{\thry{T'}}\to\topo{M}_{\theory}\]
and the induced geometric morphism
\[f:\Eqsheav{I_{\thry{T'}}}{M_{\thry{T'}}}\to \Eqsheav{I_{\thry{T}}}{M_{\thry{T}}}\]
\end{Remark}
\subsection{The theory functor \Form}
\label{Subsection: The Theory Functor}
We construct an adjoint to the contravariant semantical groupoid functor
$\Mod:\cat{T}\to\alg{Gpd}/\topo{S}$. As
noted in Section \ref{Subsection: The Semantical Groupoid Functor}, any morphism of topological groupoids
$f:\topo{G}\to\topo{S}$ induces a geometric morphism $f:\Eqsheav{G_1}{G_0}\to\Eqsheav{S_1}{S_0}$ which classifies an object $U_{\topo{G}}:=f^*(\cat{U})$ in \Eqsheav{G_1}{G_0}, with finite
powers $1,U_{\topo{G}},U^2_{\topo{G}},\ldots$ obtained by taking the (canonical) finite fiberwise products in
\Eqsheav{G_1}{G_0}. It is clear that a subobject of an equivariant sheaf can be represented (uniquely) by an open stable subset, the subset inclusion of which can be regarded as a canonical inclusion for the subobject.
\begin{Definition}
For $g:\topo{G}\to \topo{S}$ a non-empty topological groupoid over $\topo{S}$, let
\[\Form(\topo{G})\embedd \Eqsheav{G_1}{G_0}\]
be the full subcategory consisting of the (canonical) finite powers $1,U_{\topo{G}},U^2_{\topo{G}},\ldots$ of
$U_{\topo{G}}=g^*(\cat{U})$ in \Eqsheav{G_1}{G_0} together  with all equivariant sheaves $A$ such that there exists an arrow $A\to
U^{n}_{\topo{G}}$ with underlying continuous function a subspace inclusion (thus $\Form(\topo{G})$ is the full subcategory of the `relations' on $U_{\topo{G}}$).
\end{Definition}
$\Form(\topo{G})$ is then a geometric category with a distinguished object and a system of inclusions, i.e.\ an object in \cat{T}. For any morphism
\[
\bfig
\Vtriangle(1000,-300)<300,300>[\topo{H}`\topo{G}`\topo{S};f`g`h] \efig
\]
the inverse image part  $f^*:\Eqsheav{G_1}{G_0}\to\Eqsheav{H_1}{H_0}$ of the induced geometric morphism can be assumed to preserve the (finite powers of the) distinguished object as well as inclusions on the nose, and so restricts to a morphism $\Form(\topo{G})\to\Form(\topo{H})$ in \cat{T}.
\begin{Definition}\label{Definition: Form}
The contravariant functor $\Form:\alg{Gpd}/\topo{S}\to\cat{T}$ sends a non-empty groupoid \topo{G} over \topo{S} to the category $\Form(\topo{G})$ in \cat{T} and a morphism of groupoids over \topo{S}
\[
\bfig
\Vtriangle(1000,-300)<300,300>[\topo{H}`\topo{G}`\topo{S};f`g`h] \efig
\]
to the restriction to $\Form(\topo{G})$ of the inverse image $f^*:\Eqsheav{G_1}{G_0}\to\Eqsheav{H_1}{H_0}$ of the induced geometric morphism. The \emph{empty} topological groupoid,
\[\emptyset\to\topo{S}\]
is sent to the syntactic category \synt{C}{T_{\bot}} of the inconsistent geometric theory (i.e.\ the category consisting of the natural numbers with a unique isomorphism between any two of them), and the unique morphism $\emptyset\to\topo{G}$ is sent to the unique \cat{T}-morphism $\Form(\topo{G})\to \synt{C}{T_{\bot}}$.
\end{Definition}
\begin{Lemma}\label{Lemma: Form lands in TA}
For any topological groupoid $g:\topo{G}\to\topo{S}$, the category/theory $\Form(\topo{G})$ has enough \thry{S}-indexed models.
\begin{Proof} The inverse images of the points $p_x:\Sets\to\Eqsheav{G_1}{G_0}$ induced by elements $x\in G_0$ are jointly conservative and each induce an \thry{S}-indexed model of $\Form(\topo{G})$.
\end{Proof}
\end{Lemma}
Consider a  geometric theory \theory. We form the (possibly inconsistent) quotient theory $\theory_{\thry{S}}$ by adding all geometric sequents (over the signature of \theory) which are true in all models with underlying set a quotient of a subset of \thry{S}, i.e.\ true in all \thry{S}-indexed \theory-models,
\[\theory_{\thry{S}}:=\cterm{\sigma}{\alg{M}\vDash\sigma\ \textnormal{for all \thry{S}-indexed \theory-models}\ \alg{M}}\]
\begin{Definition}\label{Definition: Ta}
Define $\cat{T}_{\thry{S}}$ to be the full subcategory of \cat{T} consisting of those theories which have enough \thry{S}-indexed models, i.e.\ such that $\theory_{\thry{S}}=\theory$.
\end{Definition}
The functor $\Form$ factors through $\cat{T}_{\thry{S}}$ by Lemma \ref{Lemma: Form lands in TA}, and clearly $\Mod(\theory)=\Mod(\theory_{\thry{S}})$, i.e.\  $\topo{M}_{\theory}=\topo{M}_{\theory_{\thry{S}}}$.
Moreover,  the interpretation of a theory \theory\ into $\theory_{\thry{S}}$ yields a canonical morphism
\begin{equation}\label{Equation: Interpretation etas}
\synt{C}{T}\to^{\eta_{\theory}} \synt{C}{T_{\thry{S}}}
\end{equation}
in \cat{T} from an object in \cat{T} to an object in $\cat{T}_{\thry{S}}$, which is the unit component of an adjunction:
\begin{Lemma}\label{Lemma: T to TA adjunction}
The inclusion $\cat{T}_{\thry{S}}\to \cat{T}$ is right adjoint to the functor which sends a theory \synt{C}{T} to $\synt{C}{T_\thry{S}}$.
\begin{Proof} Straightforward.
\end{Proof}
\end{Lemma}
We therefore restrict attention to $\cat{T}_{\thry{S}}$ in order to show that there is a contravariant adjunction
\[\bfig
\morphism|b|/{@{>}@/_1em/}/<1250,0>[\cat{T}_{\thry{S}}^{op}`\alg{Grp}/\topo{S};\Mod]
\morphism/{@{<-}@/^1em/}/<1250,0>[\cat{T}_{\thry{S}}^{op}`\alg{Grp}/\topo{S};\Form]
\place(625,0)[\bot]
\efig\]
which composed with the adjunction of Lemma \ref{Lemma: Form lands in TA} yields the adjunction between \cat{T} and $\alg{Grp}/\topo{S}$. Assume in what follows, therefore, that all theories have enough \thry{S}-indexed models. Moreover, we leave the special case of the inconsistent theory implicit (recall that we defined \Form\ such that $\Form(\Mod(\synt{C}{T_{\bot}}))=\Form(\emptyset )=\synt{C}{T_{\bot}}$ for $\synt{C}{T_{\bot}}$ the inconsistent theory).
\begin{Lemma}\label{Lemma: Singlesort, f* commutes with Mdag}
The square
\[\bfig
\square<700,400>[\synt{C}{T}`\synt{C}{T'}`\Eqsheav{I_{\thry{T}}}{M_{\thry{T}}}`\Eqsheav{I_{\thry{T'}}}{M_{\thry{T'}}};F`\cat{M}_{\theory}`\cat{M}_{\thry{T'}}`f^*]
\efig\]
commutes.
\begin{Proof}
Consider, for an object $\syntob{\alg{x}}{\phi}$ in \synt{C}{T}, the square
\[\bfig
\square/<-`>`>`<-/<700,400>[\sox{\alg{x}}{\phi}`\csem{\alg{x}}{F(\phi)}_{M_{\thry{T'}}}`M_{\thry{T}}`M_{\thry{T'}};```f_0]
\efig\]
Since $f_0$ is composition with $F$, the fiber $\csem{\alg{x}}{F(\phi)}^{\alg{M}}$ over
$\alg{M}\in M_{\thry{T'}}$ is the fiber $\csem{\alg{x}}{\phi}^{f_0(\alg{M})}$ over
$f_0(\alg{M})\in M_{\theory}$, so the square is a pullback of sets. A basic open
$\bopen{\syntob{\alg{x},\alg{y}}{\psi},\alg{a}}$ is pulled back to a basic open $\bopen{\syntob{\alg{x},\alg{y}}{F(\psi)},\alg{a}}$,
so the pullback topology is contained in the logical topology. For an element \pair{\alg{M},\alg{[a]}} in basic
open $\bopen{\syntob{\alg{x},\alg{y}}{\psi},\alg{b}}\subseteq\csem{\alg{x}}{F(\phi)}_{M_{\thry{T'}}}$, the set $V=\bopen{\syntob{\alg{x},\alg{y}}{\alg{x}=\alg{y}},\alg{a}}\subseteq
\sox{\alg{x}}{\phi}$ is open and $\pair{\alg{M},\alg{[a]}}\in
V\times_{M_{\theory}}\bopen{\syntob{\alg{x},\alg{y}}{\psi},\alg{a},\alg{b}}\subseteq \bopen{\syntob{\alg{x},\alg{y}}{\psi},\alg{b}}$, so the logical topology is
contained in the pullback topology. With $f_1:I_{\thry{T'}}\rightarrow I_{\theory}$ being just a restriction
function, we conclude that $f^*\circ \cat{M}_{\theory}=\cat{M}_{\thry{T'}}\circ F$.
\end{Proof}
\end{Lemma}
\begin{Lemma} \label{Lemma: Theta after Gamma gives the same back}
For any (consistent) $\synt{C}{T}$ in $\cat{T}_{\thry{S}}$, the functor $\cat{M}:\synt{C}{T}\to \Eqsheav{I_{\theory}}{M_{\theory}}$ factors as an isomorphism (followed by an inclusion) through $\Form(\topo{M}_{\theory})$,
\[\bfig
\dtriangle/<-`^{(}->`->/<750,400>[\Form(\topo{M}_{\theory})=\Form\circ \Mod(\synt{C}{T})`\synt{C}{T}`\Eqsheav{I_{\theory}}{M_{\theory}};\epsilon_{\theory}``\cat{M}]
\dtriangle|bla|/<-`^{(}->`->/<750,400>[\Form(\topo{M}_{\theory})=\Form\circ \Mod(\synt{C}{T})`\synt{C}{T}`\Eqsheav{I_{\theory}}{M_{\theory}};\cong``]
\efig\]
\begin{Proof}
The functor $\cat{M}:\synt{C}{T}\to \Eqsheav{I_{\theory}}{M_{\theory}} \cong \mathrm{Set}[\theory]$ is an isomorphism on its image, by
Theorem \ref{Theorem: The classifying topos of T is equivariant G-sheaves} and preserves the distinguished object and inclusions
by construction. The image is $\Form(\topo{M}_{\theory})$ by an application of Lemma \ref{Lemma: Singlesort, f* commutes with Mdag}.
\end{Proof}
\end{Lemma}
The isomorphism obtained by factoring $\cat{M}$ through its image
\begin{equation}\label{Equation: Singlesort, counit component}
\bfig
\dtriangle/<-`^{ (}->`>/<750,400>[\Form(\topo{M}_{\theory})`\synt{C}{T}`\Eqsheav{I_{\theory}}{M_{\theory}};
\epsilon_{\theory}``\cat{M}_{\theory}]
\efig
\end{equation}
is our counit component candidate at \theory\ (except in the inconsistent case where it is the identity $\Form(\Mod(\synt{C}{T_{\bot}}))=\Form(\emptyset )=\synt{C}{T_{\bot}}$).
\begin{Lemma}\label{Lemma: Singlesort, Counit component candidates form nat. trans}
There is a natural transformation,
\[\epsilon:1_{\cat{T}}\to \Form\circ \Mod\]
whose component at an object \synt{C}{T} in $\cat{T}_{\thry{S}}$ is the isomorphism
\[\epsilon_{\theory}:\synt{C}{T}\to\Form(\topo{M}_{\theory})\]
of (\ref{Equation: Singlesort, counit component}).
\begin{Proof}
By Lemma \ref{Lemma: Singlesort, f* commutes with Mdag}.
\end{Proof}
\end{Lemma}

Next, we construct the unit. For a given topological groupoid $h:\topo{H}\to\topo{S}$ over
$\topo{S}$, write $\synt{C}{T}\cong \Form(\topo{H})$ just to simplify notation somewhat. We construct the object component $\eta_0:H_0\rightarrow M_{\theory}$ of  the unit
\[
\bfig \Vtriangle<350,350>[\topo{H}`\topo{M}_{\theory}`\topo{S};\eta_{\topo{H}}`h`u] \efig
\]
as follows.
%
%
%
%
%
%
%
%
%
%
%
%
%
%
%
%
Consider an element $x\in H_0$. It corresponds to a point
\[\Sets \to \Sh{H_0} \epi \Eqsheav{H_1}{H_0} \]
and thereby to a \theory-model ($\Form(\topo{H})$-model)
$\alg{M}_x$ with underlying set a quotient of \thry{S}. Specifically, the underlying set is $h_0(x)$.
Accordingly, set $\eta_0(x)=\alg{M}_x$.
%
\begin{Lemma}\label{Lemma: Singlesort, eta0 is continuous}
The map $\eta_0:H_0\rightarrow M_{\theory}$ is continuous.
\begin{Proof}
Given a basic open $V\subseteq M_{\theory}$, we may think of it as being presented in terms of an inclusion
$C\embedd U_{\topo{H}}^n$ in $\synt{C}{T}\cong \Form(\topo{H}) \embedd \Eqsheav{H_1}{H_0}$ and a tuple of elements
$\alg{a}\in\thry{S}$ with length $n$, as
\[
V=\cterm{\alg{M}\in M_{\theory}}{\alg{[a]}\in\alg{M}(C)\subseteq \alg{M}(U_{\topo{H}}^n)}.
\]
Pulled back along $\eta_0$, this is the set of those $x\in H_0$ such that there is a tuple of equivalence classes of the form
$\alg{[a]}$ in the fiber over $x$ in $C$. But this set is open, for it is the image along
the local homeomorphism $e:U^n\to H_0$ of the pullback of $W:=\bopen{\syntob{\alg{y},\alg{z}}{\alg{y}=\alg{z}}, \alg{a}}$ along
$h_0:H_0\to S_0$ intersected with $C$:
\[
\bfig
 \square<900,350>[U^n`\sox{\alg{y}}{\top}`H_0`S_0;`e``h_0]
 \place(100,250)[\pbangle]
 \square(0,350)/>`^{ (}->`^{ (}->`>/<900,350>[h_0^*(W)`W`U^n`\sox{\alg{y}}{\top};
 ```]
 \place(100,600)[\pbangle]
 \square(-900,350)/^{ (}->`^{ (}->`^{ (}->`^{ (}->/<900,350>[C\cap h_0^*
 (W)`h_0^*(W)`C`U^n;```]
\place(-800,600)[\pbangle]
\efig
\]
so that $\eta_0^{-1}(V)=\exists_{e}(C\cap h_0^*(W))$.
\end{Proof}
\end{Lemma}
Next, a point $a:x\rightarrow y$ in $H_{1}$ gives us a \theory-isomorphism between $\eta_0(x)$ and
$\eta_0(y)$---the underlying function of which is the bijection $h_1(a):h_0(x)\rightarrow h_0(y)$---and so we obtain a
function $\eta_1:H_1\rightarrow I_{\theory}$ over $S_1$,
\begin{equation}\label{Equation: Construction of unit: Commuting triangle of continuous maps of spaces of isomorphisms }
\bfig \Vtriangle<350,350>[H_1`I_{\theory}`S_1;\eta_{1}`h_1`u_{1}] \efig
\end{equation}
%
%
\begin{Lemma}\label{Lemma: Singlesort, Constructing eta}
The map $\eta_1:H_1\rightarrow I_{\theory}$ is continuous, and $\eta_1$ together with $\eta_0$ constitute a
morphism of topological groupoids $\eta_{\topo{H}}:\topo{H}\to\topo{M}_{\theory}$ over $\topo{S}$.
\begin{Proof}
By Lemma \ref{Lemma: Singlesort, eta0 is continuous} and the continuity of $h_1:H_1\rightarrow S_1$.
\end{Proof}
\end{Lemma}
\begin{Lemma}\label{Lemma: The comparison functor goes to an isomorphism}
Given a topological groupoid $h:\topo{H}\to\topo{S}$, if we apply the theory functor $\Form$ to the morphism $\eta_{\topo{H}}:\topo{H}\to \Mod\circ\Form(\topo{H})=:\topo{M}_{\theory}$, then
\[
\Form(\eta_{\topo{H}}): \Form(\topo{M}_{\theory})\to \Form(\topo{H}) \]
is left inverse to $\epsilon_{\Form(\topo{H})}:\Form(\topo{H})\to \Form(\topo{M}_{\theory})$.
%
%
%
%
\begin{Proof}
By construction.
\end{Proof}
\end{Lemma}
\begin{Lemma}\label{Lemma: EtaG is natural}
The morphism of topological groupoids over $\topo{S}$
\[
\bfig
\Vtriangle<350,350>[\topo{G}`\Mod(\Form(\topo{G}))`\topo{S};\eta_{\topo{G}}`g`u_{\Form(\topo{G})}]
\efig
\]
is natural in \topo{G}.
\begin{Proof}
Given a morphism of topological groupoids, with their induced geometric morphisms of topoi,
\[
\bfig \Vtriangle<350,350>[\topo{G}`\topo{H}`\topo{S};f`g`h]
\Vtriangle(1500,0)<350,350>[\Eqsheav{G_1}{G_0}`\Eqsheav{H_1}{H_0}`\Eqsheav{S_1}{S_0};f`g`h]\efig
\]
we need to verify that the following squares commute:
\[
\bfig
\square<750,500>[G_0`\Mod \circ \Form(\topo{G})_0` H_0` \Mod \circ \Form(\topo{H})_0; \eta_{0}`f_0`\Mod \circ
\Form(f)_0`\eta_0]
\square(1500,0)<750,500>[G_1`\Mod \circ \Form(\topo{G})_1` H_1` \Mod \circ \Form(\topo{H})_1; \eta_{1}`f_1`\Mod
\circ \Form(f)_1`\eta_1]
\efig
\]
(abusing notation somewhat). We do the left square: given an element $x\in G_0$, $\eta_0(x)$ is the
$\Form(\topo{G})$-model which sends an $C\embedd U_{\topo{G}}^n$  in $\Form(\topo{G})$ to
the fibre of $C$ over $x$. Applying $\Mod \circ \Form(f_0)$ means composing this model with
$f^*\upharpoonright_{\Form(\topo{H})}:\Form(\topo{H})\to \Form(\topo{G})$ to obtain the
$\Form(\topo{H})$-model which sends an $D\embedd U_{\topo{H}}^n$ in $\Form(\topo{H})$ to
the fibre of $f^*(D)$ over $x$. But this is precisely the fibre of $D$ over $f_0(x)$. So $\Mod \circ
\Form(f_0)\circ \eta_0=\eta_0\circ f_0$. The right square is similar.
%
\end{Proof}
\end{Lemma}
\begin{Proposition}\label{Proposition: Singlesort, Theta is left adjoint to Gamma}
$\Form$ is left adjoint to $\Mod$,
\[
\bfig
\morphism|a|/{@{<-}@/^10pt/}/<800,0>[\cat{T}_{\thry{S}}^{\mng{op}}` \alg{Gpd}/\topo{S}; \Form]
\morphism|b|/{@{>}@/^-10pt/}/<800,0>[\cat{T}_{\thry{S}}^{\mng{op}}` \alg{Gpd}/\topo{S}; \Mod] \place(375,-25)[\bot]
\efig
\]
\begin{Proof}
We need to verify the triangle identities,
\[
\bfig
\btriangle/<-`<-`<-/<1200,500>[\Form(\topo{G})` \Form\circ \Mod \circ \Form(\topo{G})` \Form(\topo{G});
\Form(\eta_{\topo{G}})` 1_{\Form(\topo{G})}` \varepsilon_{\Form(\topo{G})}]
\place(200,200)[=]
\btriangle(0,-1000)<1300,500>[\Mod(\synt{C}{T})=\topo{M}_{\theory}` \Mod\circ \Form \circ
\Mod(\synt{C}{T})` \Mod(\synt{C}{T})=\topo{M}_{\theory}; \eta_{\Mod(\synt{C}{T})}`
1_{\Mod(\synt{C}{T})}` \Mod(\varepsilon_{\Mod(\synt{C}{T}})]
\place(200,-800)[=] \efig
\]
That the top triangle commutes is Lemma \ref{Lemma: The comparison functor goes to an isomorphism}. The bottom
triangle is equally straightforward.
%
%
%
%
%
\end{Proof}
\end{Proposition}
We compose with the adjunction between \cat{T} and $\cat{T}/\thry{S}$ to obtain an adjunction between \cat{T} and  $\alg{Gpd}/\topo{S}$.
\begin{Corollary}\label{Corollary: Theta is left adjoint to Gamma}
$\Form:\alg{Gpd}/\topo{S}\to\cat{T}^{\mng{op}} $ is left adjoint to $\Mod:\cat{T}^{\mng{op}}\to\alg{Gpd}/\topo{S}$,
\[
\bfig
\morphism|a|/{@{<-}@/^10pt/}/<800,0>[\cat{T}^{\mng{op}}` \alg{Gpd}/\topo{S}; \Form]
\morphism|b|/{@{>}@/^-10pt/}/<800,0>[\cat{T}^{\mng{op}}` \alg{Gpd}/\topo{S}; \Mod] \place(375,-25)[\bot]
\efig
\]
\begin{Proof}
By Proposition \ref{Proposition: Singlesort, Theta is left adjoint to Gamma} and Lemma \ref{Lemma: T to TA adjunction}.
\end{Proof}
\end{Corollary}
\begin{Definition}\label{Definition: Singlesort, Sem}
Let \alg{Sem} be the image of $\Mod$ in $\alg{Gpd}/\topo{S}$, i.e.\ the subcategory of semantical groupoids over
$\topo{S}$.
\end{Definition}
From Lemma \ref{Lemma: Theta after Gamma gives the same back} and the remarks preceding it, it is clear that
\alg{Sem} is a full subcategory of $\alg{Gpd}/\topo{S}$. We record the following as a now easy consequence of the above:
\begin{Corollary}\label{Corollary: Gamma is an equivalence on its image}
The adjunction $\Mod \vdash \Form$ restricts to an equivalence \[ \cat{T}_{\thry{S}}^{\mng{op}}\simeq \alg{Sem}.\]
\end{Corollary}
\subsection{Characterization of semantic groupoids}
\label{Subsection: Characterization of semantic groupoids}

Corollary \ref{Corollary: Gamma is an equivalence on its image} yields a duality between the `category of theories' $\cat{T}_{\thry{S}}$ and a category \alg{Sem} of `semantic' groupoids (over \topo{S}). However, \alg{Sem} is defined as the image of a functor and lacks an independent characterization. It also seems overly restrictive to call only those groupoids in the image of $\Mod$ `semantical'. In this section, two conditions for a groupoid $f:\topo{G}\rightarrow\topo{S}$ to be a `groupoid of \thry{S}-models' are proposed. First, that it should be open and `closed under bijections', in the way a bijection $A\to\alg{M}$ into a model induces an isomorphic model on the set $A$. And second, that the topos \Eqsheav{G_1}{G_0} should classify $\Form(\topo{G})$. Using Moerdijk's site description, and in the presence of the former condition,  the latter condition can be formulated intrinsically for $f:\topo{G}\rightarrow\topo{S}$. The result is a characterization of a full subcategory $\alg{Sem}_{\thry{S}}$ of $\alg{Gpd}/\topo{S}$ which contains \alg{Sem} and which has the property that the adjunction of Proposition \ref{Proposition: Singlesort, Theta is left adjoint to Gamma} restricts to an adjunction the unit components of which are Morita equivalences of groupoids (where a morphism of topological groupoids is a Morita equivalence if the induced geometric morphism of equivariant sheaf toposes is an equivalence of categories).

Recall from Section \ref{Subsection: Equivariant sheaves on topological groupoids} the Moerdijk site of an equivariant sheaf topos. If $f:\topo{H}\rightarrow \topo{G}$ is a morphism of open topological groupoids and $N\subseteq G_1$ is open and closed under compositions and inverses, then so is $f_1^{-1}(N)\subseteq H_1$. Therefore, \pair{\topo{H},f_0^{-1}(U),f_1^{-1}(N)} is an object in the Moerdijk site of \Eqsheav{H_1}{H_0}. Moreover, it is straightforward to verify that
\[\pair{\topo{H},f_0^{-1}(U),f_1^{-1}(N)}=f^*(\pair{\topo{G},U,N})\]
\emph{if} $f$ satisfies the condition that for all $(h:x\rightarrow f_0(y))\in G_1$ there exists $g\in H_1$ such that $c(g)=y$ and $f_1(g)=h$.
%
%
Say that $f:\topo{H}\rightarrow \topo{G}$ is \emph{strongly full} if it satisfies this condition. Note that the forgetful morphism $u_{\theory}:\topo{M}_{\theory}\rightarrow \topo{S}$ which sends a model to its underlying set is strongly full, since a bijection $A\rightarrow M$ from a set into a structure/model induces an isomorphic structure on $A$.

Now, consider the groupoid \topo{S} and the classifying object $\cat{U}\in \Eqsheav{S_1}{S_0}$. For $k\geq 0$, let \alg{a} be a list of \emph{distinct} elements of \thry{S} of length $k$. Then $\bopen{\alg{a}\mapsto\alg{a}}\subseteq S_1$ is open and closed under composition and inverses, and it follows from Lemma \ref{Lemma: Definables are GUNS} that
\[\cat{U}^k\cong \pair{\topo{S},\bopen{\alg{a}},\bopen{\alg{a}\mapsto\alg{a}}}\]
in \Eqsheav{S_1}{S_0}.
\begin{Definition}\label{Definition: Groupoid of A models}
Say that a topological groupoid $f:\topo{G}\rightarrow \topo{S}$ over \topo{S} is a \emph{groupoid of \thry{S}-models} if \topo{G} is open and the following conditions are satisfied:
\begin{enumerate}[i.]
\item $f$ is strongly full; and
\item for each open $N\subseteq G_1$ which is closed under composition and inverses, and for each $x\in d(N)=c(N)$, there exists an open neighborhood $W$ of $x$ contained in $d(N)$ and a list   of distinct elements $\alg{a}\in\thry{S}$ such that $W\subseteq f_0^{-1}(\bopen{\alg{a}})$ and $f_1^{-1}(\bopen{\alg{a}\mapsto\alg{a}})\cap (d^{-1}(W)\cap c^{-1}(W))\subseteq N$.
\end{enumerate}
Let $\alg{Sem}_{\thry{S}}$ denote the full subcategory of $\alg{Gpd}/\topo{S}$ of groupoids of \thry{S}-models.
\end{Definition}
\begin{Lemma}\label{Lemma: Groupoids of models over A are semantical}
Let $f:\topo{G}\rightarrow \topo{S}$ be a groupoid of\  \thry{S}-models. Then \Eqsheav{G_1}{G_0} classifies $\Form(\topo{G})$.
\begin{Proof}Since $f$ is strongly full, we have that
\[f^*(\cat{U}^k)\cong\pair{\topo{G},f_1^{-1}(\bopen{\alg{a}\mapsto\alg{a}}),f_0^{-1}(\bopen{\alg{a}}})\]
for all lists \alg{a} of $k$ distinct elements of \thry{S}. Denote $\Form(\topo{G})$ by \theory\ and recall that it consists of the subcategory spanned by subobjects of $f^*(\cat{U}^k)$, for $k\geq 0$. By Lemma \ref{Lemma: GUN subobjects}, the Moerdijk site is closed under subobjects. Thus \theory\ is, up to isomorphism, a full subcategory of the Moerdijk site of \Eqsheav{G_1}{G_0}. We show that \theory\ is generating. Let \pair{\topo{G},U,N} be an object in the Moerdijk site of  \Eqsheav{G_1}{G_0}, and let $[f:x\rightarrow y]$ be an element of  \pair{\topo{G},U,N}. Chose an open set $W$ and a list \alg{a} of distinct elements of \thry{S} such that $x\in W\subseteq U$;  $W\subseteq f_0^{-1}(\bopen{\alg{a}})$; and $f_1^{-1}(\bopen{\alg{a}\mapsto\alg{a}})\cap (d^{-1}(W)\cap c^{-1}(W))\subseteq N$. Set $M:= f_1^{-1}(\bopen{\alg{a}\mapsto\alg{a}})\cap (d^{-1}(W)\cap c^{-1}(W))$ and consider $\pair{\topo{G},W,M}$. First, the obvious morphism
\[\pair{\topo{G},W,M}\rightarrow \pair{\topo{G},f_0^{-1}(\bopen{\alg{a}}),f_1^{-1}(\bopen{\alg{a}\mapsto\alg{a}})}\cong f^*(\cat{U}^k)\]
is injective, so \pair{\topo{G},W,M} is in \theory. Second, sending an equivalence class $[h]_{\sim_M}$ to $[h]_{\sim_N}$ well-defines a morphism
\[\pair{\topo{G},W,M}\rightarrow \pair{\topo{G},U,N}\]
with $[f:x\rightarrow y]$ in its image. So restricting the Moerdijk site of \Eqsheav{G_1}{G_0} to the full subcategory \theory\ still yields a generating set.
\end{Proof}
\end{Lemma}
\begin{Corollary}\label{Corollary: Sem+ Morita equivalent units}
Let $f:\topo{G}\rightarrow \topo{S}$ be a groupoid of \thry{S}-models. Then the unit $\topo{G}\rightarrow \Mod(\Form(\topo{G}))$ is a Morita equivalence.
\end{Corollary}
\begin{Lemma}\label{Lemma: Mod factros through Sem+}
The functor $\Mod:\cat{T}^{\mng{op}}\to\alg{Gpd}/\topo{S}$ factors through $\alg{Sem}_{\thry{S}}$.
\begin{Proof} By Proposition \ref{Lemma: The definables generate equivariant sheaves on the groupoid of models}.
\end{Proof}
\end{Lemma}
\begin{Theorem}
The adjunction $\Mod \vdash \Form$ restricts to an adjunction
\[
\bfig
\morphism|a|/{@{<-}@/^10pt/}/<800,0>[\cat{T}^{\mng{op}}_{\thry{S}}` \alg{Sem}_{\thry{S}}; \Form]
\morphism|b|/{@{>}@/^-10pt/}/<800,0>[\cat{T}^{\mng{op}}_{\thry{S}}` \alg{Sem}_{\thry{S}}; \Mod] \place(375,-25)[\bot]
\efig
\]
the counit and unit components of which are isomorphisms and Morita equivalences, respectively.
\end{Theorem}
\subsection{Coherent theories}\label{Subsection: Coherent Theories}
We end with a compact note on \emph{coherent} theories, that is, theories that can be axiomatized using only sequents involving finitary formulas (no infinite disjunctions). A topos is coherent if it classifies a coherent theory, and a category is coherent if it is regular and has finite, stable unions of subobjects (see \cite{elephant1} for more on coherent theories, coherent categories, and coherent toposes). All coherent theories have enough models, although not, of course, necessarily enough \thry{S}-indexed models for a fixed index set \thry{S}. We represent coherent theories in the same way as geometric theories were represented in Section \ref{Subsection: Theories and Categories}, restricting to single-sorted coherent theories that, for the sake of brevity, have enough \thry{S}-indexed models (for some fixed index set \thry{S}). For multi-sorted theories, see instead the approach in \cite{awodeyforssell}.   Let $\cat{T}^c_{\thry{S}}$ be the category of (small) coherent categories with enough \thry{S}-indexed models and with a distinguished object and a system of inclusions, and with distinguished object-preserving coherent functors between them. There is an obvious functor
\[\cat{T}^c_{\thry{S}}\to \cat{T}\]
the image of which are the coherent theories with enough \thry{S}-indexed models. Composing with \Mod, we obtain a functor
\[\Mod^c :(\cat{T}^c_{\thry{S}})^{\mathrm{op}}\to\alg{Sem}_{\thry{S}}\]
Say that a frame is coherent if it is compact; compact elements are closed under meets; and every element is a join of compact elements (as in \cite{johnstone:82}). If $f:\topo{G}\rightarrow\topo{S}$ is a  groupoid of \thry{S}-models in the image of $\Mod^c(\cat{T}^c_{\thry{S}})$ then (i) the frame of subobjects of $\cat{U}^k_{\topo{G}}$ is coherent for $k\geq 0$ and (ii)  pullbacks of  compact subobjects along  projection maps $\cat{U}^{k+1}_{\topo{G}}\rightarrow \cat{U}^k_{\topo{G}}$ are again compact. Conversely, these two conditions (entailing that the category of compact relations on $\cat{U}_{\topo{G}}$ is generating and coherent) are sufficient for  $\Form(\topo{G})$ to be coherent for a groupoid of \thry{S}-models $f:\topo{G}\rightarrow\topo{S}$. Using Moerdijk's site description, conditions (i) and (ii) can be expressed in terms of the groupoids \topo{G} and \topo{S} and the morphism $f$. First, since $f$ is strongly full, we have
\[\cat{U}^k_{\topo{G}}\cong \pair{\topo{G},f_0^{-1}(\bopen{\alg{a}}),f_1^{-1}(\bopen{\alg{a}\mapsto\alg{a}})}\]
for some (any) list of $k$ distinct elements \alg{a} in \thry{S}. So the frame of subobjects of $\cat{U}^k_{\topo{G}}$ is coherent iff the frame of open subsets of $f_0^{-1}(\bopen{\alg{a}})$ that are closed under $f_1^{-1}(\bopen{\alg{a}\mapsto\alg{a}})$ is  coherent. Second  it is straightforward to compute (see \cite{forssell:subgpds}) that if \alg{a} is a list of $k+1$ distinct elements and \alg{b} is a sublist of $k$ of them, then the projection map
\[\cat{U}^{k+1}\cong \pair{\topo{S},\bopen{\alg{a}},\bopen{\alg{a}\mapsto\alg{a}}}\to \cat{U}^{k}\cong \pair{\topo{S},\bopen{\alg{b}},\bopen{\alg{b}\mapsto\alg{b}}} \]
can be given as precomposing with arrows from the open set
\[T=\left(\begin{array}{c}
\bopen{\alg{b}}  \\
  \alg{b} \mapsto \alg{b}   \\
 \bopen{\alg{a}}
\end{array}\right)\]
and that the projection map
\[\cat{U}^{k+1}_{\topo{G}}\cong \pair{\topo{G},f_0^{-1}(\bopen{\alg{a}}),f_1^{-1}(\bopen{\alg{a}\mapsto\alg{a}})}\to<100> \cat{U}^k_{\topo{G}}\cong \pair{\topo{G},f_0^{-1}(\bopen{\alg{b}}),f_1^{-1}(\bopen{\alg{b}\mapsto\alg{b}})}\]
can similarly be obtained by precomposing with arrows from $f_1^{-1}(T)$. Using Lemma \ref{Lemma: GUN subobjects}, condition (ii) then spells out as: for any open subset $S\subseteq G_0$  that is a compact element of  the frame of open subsets of $f_0^{-1}(\bopen{\alg{b}})$ that are closed under $f_1^{-1}(\bopen{\alg{b}\mapsto\alg{b}})$, the open subset
\[\cterm{x\in f_0^{-1}(\bopen{\alg{a}})}{\fins{h\in f_1^{-1}(T)}c(h)=x\wedge d(h)\in S}\]
is compact in the frame of open subsets of $f_0^{-1}(\bopen{\alg{a}})$ that are closed under $f_1^{-1}(\bopen{\alg{a}\mapsto\alg{a}})$.
\begin{Definition}\label{Definition: Coherent groupoids of models}
Say that $f:\topo{G}\rightarrow\topo{S}$ is a \emph{coherent groupoid of \thry{S}-models} if it is a groupoid of \thry{S}-models such that conditions (i) and (ii) above are satified. Let
\[\alg{CSem}_{\thry{S}}\embedd \alg{Sem}_{\thry{S}}\]
be the full subcategory of coherent groupoids of \thry{S}-models.
\end{Definition}
Given a coherent groupoid of \thry{S}-models $f:\topo{G}\rightarrow\topo{S}$  we obtain from \Form\ a functor $\Form^c:\alg{CSem}_{\thry{S}} \to (\cat{T}^c_{\thry{S}})^{\mng{op}}$ by taking the full subcategory of compact relations on $\cat{U}_{\topo{G}}$ (instead of all relations), and we conclude:
\begin{Theorem}
The adjunction $\Mod \vdash \Form$ induces an adjunction
\[
\bfig
\morphism|a|/{@{<-}@/^10pt/}/<800,0>[(\cat{T}^c_{\thry{S}})^{\mng{op}}` \alg{CSem}_{\thry{S}}; \Form^c]
\morphism|b|/{@{>}@/^-10pt/}/<800,0>[(\cat{T}^c_{\thry{S}})^{\mng{op}}` \alg{CSem}_{\thry{S}}; \Mod^c] \place(375,-25)[\bot]
\efig
\]
between coherent theories with enough \thry{S}-indexed models and coherent groupoids of \thry{S}-models, with the counit and unit components of the adjunction being isomorphisms and Morita equivalences, respectively.
\end{Theorem}

\section*{Acknowledgements}
Thanks to Steve Awodey who has been of much help in the writing of this paper and who supervised the PhD work on which it is based. Thanks also to Ji\v{r}\'{i} Rosick\'{y} and to the Eduard \v{C}ech Center for their support through grant no.\ LC505.

\bibliographystyle{ieeetr}
\bibliography{bibliografi}
\end{document}